\documentclass[12pt,a4paper]{article}

\usepackage[utf8]{inputenc}

\usepackage[T1]{fontenc}
\usepackage{verbatim}
\usepackage[letterpaper,top=25mm,bottom=35mm,left=25mm,right=25mm]{geometry}
\usepackage{color}
\usepackage{enumerate}
\usepackage{amstext}
\usepackage{amsmath}
\usepackage{amsfonts}
\usepackage{mathrsfs}
\usepackage{amssymb, amsmath, amsthm}
\usepackage{latexsym}
\usepackage{xspace}

\usepackage{url}
\usepackage{hyperref}
\usepackage[T1]{fontenc} % pacote para hifenizar corretamente
\usepackage[all]{xy}
\usepackage{fancyhdr}
\usepackage{xspace}

\usepackage{indentfirst,amsfonts,amsmath,amsbsy,amsthm,amssymb,amscd}
\usepackage{enumerate}
\usepackage[normalem]{ulem}  
\usepackage{bbm}
\usepackage{mathrsfs}
\usepackage{dsfont}
\usepackage{hyperref}
\usepackage{makeidx}
\usepackage{calligra}
\usepackage{ulem}

\usepackage{graphicx}

\def\XXint#1#2#3{{\setbox0=\hbox{$#1{#2#3}{\int}$ }
\vcenter{\hbox{$#2#3$ }}\kern-.6\wd0}}

\newcommand{\R}{\mathbb{R}}
\newcommand{\N}{\mathbb{N}}

\newcommand{\HH}{\mathbb{H}}
\renewcommand{\H}{\mathcal{H}}
\newcommand{\taureg}{\tau_A^\mathrm{reg}}
\newcommand{\Areg}{A_0^\mathrm{reg}}
\newcommand{\lsc}{l.s.c.\ }
\newcommand{\usc}{u.s.c.\ }
	\newcommand{\ie}{\textit{i.e.}\ }
\newcommand{\HA}{$\mathrm{(HA)}$\xspace}
\newcommand{\HB}{$\mathrm{(HB)}$\xspace}
\newcommand{\HC}{$\mathrm{(HC)}$\xspace}
\newcommand{\HD}{$\mathrm{(HD)}$\xspace}

\newtheorem{theorem}{Theorem}[section]
\newtheorem{lemma}[theorem]{Lemma}
\newtheorem{remark}[theorem]{Remark}
\newtheorem{corollary}[theorem]{Corollary}
\newtheorem{proposition}[theorem]{Proposition}
\newtheorem{definition}[theorem]{Definition}

\numberwithin{equation}{section}

\newcommand{\liminfs}{\mathop{\liminf{\!}_*\,}}
\newcommand{\limsups}{\mathop{\limsup{\!}^*\,}}
\newcommand{\e}{\varepsilon}

\newcommand{\LOCa}{$(\mathrm{LOC 1})$\xspace}
\newcommand{\LOCb}{$(\mathrm{LOC 2})$\xspace}
\providecommand{\keywords}[1]
{
  \small	
  \textbf{\textit{Keywords---}} #1
}

\title{Unbounded Hamilton-Jacobi-Bellman Equations with one
co-dimensional discontinuities}
\author{Emmanuel Chasseigne$^1$, Robson Carlos Reis$^2$, Silvia Sastre-Gomez$^3$ \\
        \small $^{1}$Universit\'e de Tours, France \\
        \small $^{2}$Universidade Federal do Rio Grande do Norte, Brazil \\
        \small $^{3}$Universidad de Sevilla, Spain \\
}
\date{} % Comment this line to show today's date

\begin{document}
\maketitle

\begin{abstract}
    The aim of this work is to deal with a discontinuous Hamilton-Jacobi equation in the whole
    euclidian $N$-dimensional space, associated to a possibly unbounded optimal control problem.
    Here, the discontinuities are located on a hyperplane and the typical questions we address
    concern the existence and uniqueness of solutions, and of course the definition itself of
    solution. We consider viscosity solutions in the sense of Ishii. The convex Hamiltonians are
    associated to a control problem with specific cost and dynamics given on each side
    of the hyperplane. We assume that those are Lipschitz continuous but the main difficulty we deal
    with is that they are potentially unbounded, as well as the control spaces. Using Bellman's
    approach we construct two value functions which turn out to be the minimal and maximal solutions
    in the sense of Ishii. Moreover, we also build a whole family of value functions, which are
    still solutions in the sense of Ishii and  connect continuously the minimal solution to the
    maximal one.  
\end{abstract}
\hspace{10pt}


%==================================================%

%%%%%%%%%%%%%%%%%%%%%%%%%%%%%%%
\section{Introduction}
%%%%%%%%%%%%%%%%%%%%%%%%%%%%%%%

Inspired by the works of Kruzhkov, \textit{Viscosity Solutions'} theory has been developed and
applied to a wide extent since its birth in the early 80's. We refer the reader to
\cite{CrandalEvansLions84, CrandalLions83} for the initial papers on the subject, the famous
\textit{User's guide} \cite{CrandalIshiiLions92} and the book of Barles \cite{guy} for a more
complete overview of the theory. 

This theory is well-established and now quite stable in the case of continuous Hamiltonians,
but dealing with discontinuities has always been difficult due to the pointwise nature of the concept
of solution. And apart from several specific cases, even dealing with simple discontinuous problems
has not been fully addressed until recently. We refer to the book \cite{Barles2018} for a widespread 
introduction, information and results about discontinuous Hamilton-Jacobi problems.

In \cite{manoel} G. Barles, A. Briani and E. Chasseigne study discontinuous
Hamilton-Jacobi-Bellman equations in $\R^N$, associated to an infinite horizon control problem where
the discontinuities are located on the hyperplane $\mathcal{H}=\{x_N=0\}$. A more global approach
can be found in \cite{BarlesChasseigne2015,Barles2018} but here we
focus on this hyperplane, stationary situation.

Before presenting the approach in \cite{manoel}, let us already mention that so far,
the literature on discontinuous Hamilton-Jacobi-Bellman equations only focuses on bounded control
sets, bounded dynamics and costs. One of the main goals of the present paper is to bridge the gap to
the unbounded case, using some specific hypotheses.

\

\noindent\textsc{The bounded case} ---
In order to give a quick overview of the results in \cite{manoel} and introduce some notations and
concepts that we use throughout this work, let us decompose the space into three parts
$\mathbb{R}^{N}=\Omega_1 \cup \mathcal{H} \cup \Omega_2$, where $ \Omega_1 = \{x \in \mathbb{R}^N \
| \ x_N > 0\}$ and  $\Omega_2= \{x \in \mathbb{R}^N  \ | \ x_N < 0\}$, $\H$ being the hyperplane
separating $\Omega_1$ and $\Omega_2$: $\mathcal{H}=\{x \in \mathbb{R}^N \ | \ x_N = 0\}$.
We will take $A_1,A_2\subset\R^{d}$ to be the sets of controls. In each domain $\Omega_i$, a control
problem is defined through a dynamic function $b_i:\overline{\Omega_i}\times A_i\to\mathbb{R}^N$ and
a cost function $l_i:\overline{\Omega_i}\times A_i\to\mathbb{R}$.  As in \cite{manoel}, for the
moment the reader may assume that each $A_i$ is compact and the $(b_i,l_i)$ are at least continuous
and bounded on $\overline{\Omega}_i\times A_i$.

In order to define a value function, it is first necessary to define trajectories that may cross or
stay on $\mathcal{H}$, hence we also need to define the dynamics and cost on $\mathcal{H}$.
Following \cite{manoel}, we set $A:=A_1\times A_2 \times [0,1]$ and use the control set formed by
bounded measurable functions $\mathcal{A}:=L^{\infty}(0, \infty; A)$. Then on $\mathcal{H}$ we 
introduce a relaxed dynamic $b_\mathcal{H}:\mathcal{H}\times A\to\mathbb{R}^N$ given by a convex
combination of $b_1$ and $b_2$ through a parameter $\mu\in[0,1]$, and similarly a cost function
$l_\mathcal{H}:\mathcal{H}\times A\to\mathbb{R}$ given by a convex combination of $l_1$ and $l_2$.
More precisely, $b_\mathcal{H}$ and $l_\mathcal{H}$ are defined as 
\[
\begin{array}{ll}
	b_\mathcal{H}(x,a)=\mu b_1(x,\alpha_1) + (1-\mu) b_2(x,\alpha_2), & x\in \mathcal{H}, 
    ~a=(\alpha_1,\alpha_2,\mu)\in A,\\
	l_\mathcal{H}(x,a)=\mu l_1(x,\alpha_1) + (1-\mu) l_2(x,\alpha_2), & x\in \mathcal{H}, 
    ~a=(\alpha_1,\alpha_2,\mu)\in A.
\end{array}
\]

Following this construction, a global formulation of trajectories can be obtained by solving a
differential inclusion, we detail this approach in the Preliminaries section below.
As shown in \cite{manoel}, any trajectory $X(\cdot)$ solving the differential inclusion can be associated to a
control function $a=(\alpha_1, \alpha_2, \mu) \in \mathcal{A}:=  L^{\infty}(0,\infty;A)$ such that
\begin{equation}\label{trajec0} 
   \begin{split}
       \dot X\left(  t \right) = &  ~b_1\big(X(t), \alpha_{1}(t) \big)
       \mathds{1}_{ \lbrace X(t)\in 
       \Omega_1\rbrace} + b_2\big(X(t), \alpha_{2}t)\big) 
       \mathds{1}_{ \lbrace X_{x}(t)\in  \Omega_2\rbrace} 
    \\
    & + b_{\mathcal{H}}\big( X(t), a(t) \big)
       \mathds{1}_{ \lbrace X_{x}(t)\in 
    \mathcal{H}\rbrace} \ \ \ \ \mbox{for  a.e. } \ \ t>0\,.
    \end{split}
\end{equation}

Denoting by $\tau_A(x)$ the set of such controlled trajectories, $(X,a)$, starting from $X(0)=x$, it is
natural to introduce the following value function 
\begin{equation}\label{umenos0}
    U^{-}_{A}(x) := \inf_{(X,a)\in\tau_{A}(x)}\left( \int_{0}^{\infty} l(X(t),a(t)) 
    e^{-\lambda t} dt \right)\
\end{equation}
where $l$ stands for the global cost: $l=l_1,l_2$ or $l_\H$ depending on the location.
Introducing in each $\Omega_i$ the Hamiltonian
$$ 
    H_i(x, u, p) =\sup\limits_{\alpha_i\in A_i} \big\{ \lambda u -b_i(x,\alpha_i)\cdot p -l_i(x, \alpha_i)
    \big\}\,,
$$
it is clear enough that the value function satisfies $H_i(x,u,Du)=0$ in each $\Omega_i$---in the viscosity
sense---, a fact which can be actually obtained by performing the standard analysis through the dynamic
programming principle locally inside each $\Omega_i$.

It remains to understand the situation on $\mathcal{H}$ and, as the viscosity theory suggests, it
turns out that the value function $U^-_A$ satisfies the Ishii conditions on $\mathcal{H}$, leading
to a full set of (in)equations:
\begin{equation}\label{ISHII}  \left\{\begin{matrix} 
H_1\left(x,u, Du \right) = 0 & \mbox{in }   \Omega_1\,, 
\\ H_2\left(x, u, Du \right) = 0 & \mbox{in } \Omega_2 \,,
\\ \min\{H_1\left(x, u, Du \right), H_2\left(x, u, Du \right)\} \leq 0   & \mbox{on }
    \mathcal{H}\,,
\\ \max\{H_1\left(x, u, Du \right), H_2\left(x, u, Du \right)\}\, \geq \, 0   & \mbox{on }
    \mathcal{H}\,.
 \end{matrix} \right.
 \end{equation}
Notice that on $\H$, only the $\min$-inequality is required for a subsolution in the sense of
Ishii, while only the $\max$-one is required for a supersolution.
However, the special solution $U^-_A$ satisfies a complemented $(N-1)$-dimensional inequation on
$\mathcal{H}$: $H_T(x,u,Du)\leq0$, where the tangential Hamiltonian is defined by
\begin{equation}\label{ishiit}
    H_{T}(x,\phi, D_{\mathcal{H}}\phi)=\sup\limits_{(\alpha_1, \alpha_2,\mu)\in A_0(x)} \Big\{
        \lambda\phi(x) - b_{\mathcal{H}}(x, \alpha_1,\alpha_2, \mu)\cdot  D_{\mathcal{H}}\phi
        (x) -l_{\mathcal{H}}( x,\alpha_1, \alpha_2,\mu) \Big\}\,.
\end{equation}
Here, $A_0(x)$ is the set of controls that allow the trajectory to remain on $\mathcal{H}$, that is,
the controls such that $b_{\mathcal{H}}(x,a)\cdot e_N=0$ and
$$
 D_{\mathcal{H}}\phi \left( x\right)=\left({\partial \phi\over \partial x_1},
 \dots, {\partial \phi\over \partial x_{n-1}}\right)\,.
$$
Adding the tangential subsolution condition $H_T(x,u,Du)\leq0$ to problem \eqref{ISHII} allows to
prove a comparison result between sub and supersolutions, leading to the fact that $U^-_A$ is the
unique Ishii solution satisfying $H_T\leq0$. 

Including the subsolution condition $H_T\leq0$ leads to the notion of \textbf{stratified
solutions} which is the good notion of solution for treating \eqref{ISHII}, meaning that
a complete comparison result between sub and supersolutions holds in this framework. 

\

\noindent\textsc{Regular dynamics, extremal Ishii solutions} ---
The precise analysis of the situation performed in \cite{manoel} shows that at least two specific
value functions can be built, the ``natural'' one being of course $U^-_A$. But the authors also
build a second solution in the sense of Ishii by introducing the \textbf{regular dynamics} on
$\mathcal{H}$, which are defined as $b_\H=\mu b_1+(1-\mu)b_2$ satisfying 
$$ 
    \forall x\in\H\,,\quad b_1(x,\alpha_1)\cdot e_N \leq 0
    \quad\mbox{and}\quad b_2(x,\alpha_2)\cdot e_N \geq 0\,.
$$ 
Intuitively, such dynamics are maintaining the trajectory on $\mathcal{H}$ by ``pushing'' from
either side of the hyperplane, contrary to singular dynamics which corresponds to an equilibrium
obtained by ``pulling'' from both sides.

Defining the set of {\bf regular trajectories} as
\begin{equation}\label{controlregular} \taureg(x):=\Big\{(X(\cdot),a(\cdot)) \in 
    \tau_A(x)\ \mbox{and for a.e.} \ t\in \mathcal{E}_{\mathcal{H}},\ b_{\mathcal{H}}(X(t),a(t))
    \  \mbox{is regular}\Big\}
\end{equation}
where $\mathcal{E}_{\mathcal{H}}:=\{t: \  X(t) \in \mathcal{H} \}$, allows to construct a second
value function: 
\begin{equation}\label{umais0}
    U_{A}^{+}(x): = \inf_{(X,a)\in\taureg(x)}\left( \int_{0}^{\infty} l( X(t),a(t)) e^{-\lambda t}
    dt \right)\,.
\end{equation}
As shown in \cite{manoel}, it turns out that $U^{-}_A$ is the minimal (super)solution of the Ishii
problem \eqref{ISHII}, while $U^{+}_A$ is the maximal (sub)solution. Of course, $U^-_A$ is the
unique stratified solution.

Let us mention finally that in \cite{manoel2014} the authors extend the results in several
directions by considering general domains and finite horizon control problems under weaker
controllability assumptions. An even more global approach is performed in \cite{Barles2018} but the
authors always consider mainly bounded control sets, bounded dynamics and costs functions, apart
from some remarks on some unbounded cases like KPP-type problems.

\

\noindent\textsc{Main results of the paper} ---
As was said above, in \cite{manoel,manoel2014,Barles2018} the Ishii problem (\ref{ISHII}) is treated in
the context of compact control space and bounded cost-dynamics, restricting the study to sublinear
Hamiltonians. In this work one of the main objectives is to verify which results of \cite{manoel},
\cite{manoel2014} are still valid when dealing with non-compact controls spaces and unbounded
cost-dynamics.  

The clear motivation for considering non-compact controls spaces is  that some important
Hamiltonians can only be approached in this context, like the following typical quadratic example 
\begin{equation}\label{ex:quadratic}
    H_i(x,u,Du)=\lambda u +c_i|Du|^{2} -f_i(x)\,,\quad |f_i(x)|\leq C|x|^2\,.
\end{equation}
While the literature related to superlinear Hamiltonians is extensive --- see for example
\cite{unbounded}, \cite{representacao}, \cite{motta} \cite{unbounded2000} for control related
topics ---, dealing with {\bf discontinuities} in non-compact settings is far from being easy.

So, in this work we consider unbounded sets of controls, typically
$A_i=\mathbb{R}^d$, and only locally bounded cost-dynamic functions $(b_i,l_i)$. 
In this framework, a central hypothesis in our paper is the following:
\begin{equation}\label{alphainfinite0oo} 
    \lim_{|\alpha_i| \rightarrow \infty} \frac{l_i(y,\alpha_i)}{1+|b_i(y,\alpha_i)|}=\infty, \ \  
    \mbox{locally uniformly with respect to} \ \ y.
\end{equation}
The key idea here is that the cost associated to large dynamics is so big that such
strategies are not optimal. This allows to recover some compactness of the trajectories and the
associated controls. This type of hypothesis already appears in \cite{unbounded} where the reader
will find counter-examples to uniqueness when it is not satisfied. Using \eqref{alphainfinite0oo}
allows to reduce several arguments to the case of bounded control sets where we can use the results
of \cite{manoel} even if not all the arguments are so easy to handle.
The main results of the paper can be summarized as:

\noindent\textit{Assuming \eqref{alphainfinite0oo}, a global comparison result holds for \eqref{ISHII}
        between stratified sub and supersolutions. As a consequence, $U_A^-$ is the unique stratified
        solution of this problem.}

        \

We also build a whole family of value functions, which are still solutions in the sense of Ishii and
connect continuously the minimal solution to the maximal one. We have already commented that
problems with discontinuities, in general, do not have uniqueness. In fact, $U^{+}_A$ and $U^{-}_A$ are
viscosity solutions of the Ishii problem, but they are not the only ones. We build a family, $U^{\eta}_A$, 
of locally Lipschitz Ishii solutions. Under appropriate assumptions we obtain
that the limit  when $\eta$ goes to zero is $U_A^{+}$ and when $\eta$ goes to infinite is
$U_A^{-}$. This family can then be seen as a continuous path connecting $U^-_A$ to $U^+_A$. Such
solutions are built on the relaxed regular trajectories
\begin{equation}\label{-+controleta0} 
    \begin{aligned}\tau^{\eta}_{A}(x):=\Big\{\big(X(\cdot),a(\cdot)\big) \in \tau_A(x) &\quad
        \mbox{such that for a.e.} 
    \ t\in \mathcal{E}_{\mathcal{H}}\;,\\
        & b_2(X(t),a(t))\cdot e_N \geq -\eta,\ b_1(X(t),a(t))\cdot e_N \leq \eta \Big\}\,.
    \end{aligned}
\end{equation}
These $\eta-$trajectories may not be regular, but they are almost regular if $\eta$ is close to 0.
On the contrary, if $\eta$ is close to $+\infty$, we recover most of the trajectories. The
$\eta$-value function is then defined as one can expect: 
\begin{equation}\label{ueta0}
    U^{\eta}_{A}(x): = \inf_{(X,a)\in\tau^{\eta}_{A}(x)}\left( \int_{0}^{\infty}
    l( X(t),a(t)) e^{-\lambda t} dt \right).
\end{equation}

\

%{\color{red} Check Organization when paper stabilized}
%
%\noindent\textsc{Organization ---} In section 2 we  verify that value functions satisfy the Dynamic
%Programming Principle and are locally Lipschitz. Furthermore, using \ref{alphainfinite0oo} we show
%that the supremum defining $H_i$ is locally achieved in a compact set of controls. Consequently our
%Hamiltonians $H_i$ are well defined and are continuous functions. Moreover taking a compact sequence
%$(A_m)_m$ and considering $U_{A^{m}}^{+}$, $U_{A^{m}}^{-}$ we prove that both are solutions of an
%approximate Ishii problem. Thus, due to Half relaxed limit method sending $m$ to infinity, we
%obtain that $U_A^{+}$ and $U_A^{-}$ are viscosity solutions of the initial Ishii problem. In
%section 3, we provide growth conditions that allow us to obtain comparison and uniqueness results.
%Moreover we build a family of examples satisfying our growth conditions that guarantee uniqueness
%and we consider the Filippov approximation. In section 4, we work with value functions $U^{\eta}$.

\

\tableofcontents

%%%%%%%%%%%%%%%%%%%%%%%%%%%%%%%%%%%%%%%%%%%%%%
\section{Preliminaries}
%%%%%%%%%%%%%%%%%%%%%%%%%%%%%%%%%%%%%%%%%%%%%%

Let us begin with some basic notation, already introduced in the introduction:
throughout this paper, $A_1,A_2\subset\R^d$ are (non compact) metric spaces which are both
closed and convex, and $\Omega_1,\Omega_2,\H \subset\R^N$ are defined by
$$ \Omega_1 = \{x \in \mathbb{R}^N : x_N > 0\}\,,\quad 
\Omega_2= \{x \in \mathbb{R}^N: x_N < 0\}\,,\quad 
\H = \{x \in \mathbb{R}^N : x_N =0\}\,.$$
We assume that for $i=1,2$ we are given a pair of dynamic-cost functions $b_i,l_i$ defined on
$\overline{\Omega}_i\times A_i$ and given a fixed parameter $\lambda>0$ (the actualization factor),
we define for $x\in\Omega_i$ the Hamiltonian 
$$H_i(x,u,p):=\sup_{\alpha_i\in A_i}\big\{\lambda u-b(x,\alpha_i)\cdot p-l_i(x,\alpha_i)\big\}\,.$$
More precise assumptions on $(b_i,l_i)$ are given below.

%---------------------------------------------
\subsection{Main assumptions}

Many arguments in this paper are based on compact approximations of the control sets $A_i$. To this
end, let us introduce a sequence of compact sets $A_i^1\subset A_i^2\subset\dots\subset A_i^m\dots
\subset A_i$ such that 
\begin{equation}\label{seqqcontrol} 
    A_i= \bigcup_{m\in \mathbb{N}} A_i^m\,,
\end{equation}
We consider also all the definitions related to the control sets $A_i^{m}$: the associated
Hamiltonians $H_1^m,H_2^m$, the Ishii problem, as well as the value functions $U_{A^{m}}^{+}$,
$U_{A^{m}}^{-}$ and the tangential Hamiltonian $H_T^{m}$ which will be defined below.

Let us now state the main hypotheses on the dynamic and cost functions that we use.
\begin{enumerate}

    \item[\HA] for $i=1,2$, $(x,\alpha_i)\mapsto l_i(x,\alpha_i)$ and $(x,\alpha_i)\mapsto
        b_i(x,\alpha_i)$
        are continuous functions with respect to $(x,\alpha_i)$; $b_i$ is Lipschitz continuous with
         respect to $\alpha_i$, uniformly for $x$ in compact subsets of~$\R^N$\,; 
    \item[\HB] for $i=1,2$, any $x\in\R^N$ and $m\in\N_*$, the set 
        $\{(b_i(x,\alpha_i), l_i(x,\alpha_i)):  \alpha_i \in A^m_i \}$ is closed and convex. 
        Moreover, there exists $\delta>0$ such that
        \begin{equation*}
            \overline{B_{\delta}(0)} \subset \mathcal{B}^{A_i^m}(x), \quad\mbox{ where }
            \mathcal{B}^{A_i^m}(x):=\{b_i(x,\alpha_i): \alpha_i\in A_i^m\}\,.
        \end{equation*}
    \item[\HC] for $i=1,2$,  $l_i:\R^N\times A_i\to \R$ and $b_i:\R^N\times A_i\to \R^N$ satisfy 
        $$ \lim_{|\alpha_i|\rightarrow \infty}|b_i(x,\alpha_i)|=+\infty \ \ \mbox{and} \ \
        \lim_{|\alpha_i|\rightarrow \infty} \frac{l_i(x,\alpha_i)}{1+|b_i(x,\alpha_i)|}=+\infty\,,$$
    uniformly with respect to $x$ in compact subsets of $\R^N$\,;
\end{enumerate}

\

Hypothesis \HA is quite natural in this non-compact setting. Notice that the total controllability
assumption \HB is also automatically satisfied by $A_i$, $i=1,2$. More important is hypothesis
$\mbox{\HC}$ which states that if the dynamic grows too much, the associated
cost gets very high. Hence, the associated trajectories with high cost are not the ones that are
important in the definition of the value function, or the dynamic programming principle.

%------------------------------------------------
\subsection{The control problem}

Following \cite{manoel}, global trajectories of the control problem are defined by solving the
following differential inclusion 
\begin{equation}\label{t-bis}
    \begin{cases}
        \dot X(t)\in \mathcal{B}\big(X(t)\big)  \mbox{ for a.e.} \  t \in (0,\infty)\,,\\
        X(0) = x\, \mbox{ for } x\in\R^N,
    \end{cases}
\end{equation}
where the dynamic set-valued map is defined by
\begin{equation}\label{tr0}
 \mathcal{B}(x):= \left\{\begin{matrix} 
     \mathcal{B}_{1}(x):= \big\{b_1(x,\alpha_1): \alpha_1 \in A_1\big\}\,, & \mbox{if } x_N>0\,,\\
     \mathcal{B}_{2}(x):= \big\{b_2(x,\alpha_2): \alpha_2 \in A_2\big\}\,, & \mbox{if } x_N<0\,,\\
     \overline{\rm co}\big(\mathcal{B}_{1}(x) \cup  \mathcal{B}_{2}(x)\big)   & \mbox{if } x_N=0\,,
 \end{matrix} \right.
\end{equation}
$\overline{\rm co}(\cdot)$ being the convex hull. We then say that  $X(\cdot)$ is a {\bf
trajectory} if it is a Lipschitz continuous function that satisfies the following differential
inclusion 
\begin{equation}\label{t}
    \dot X(t)\in \mathcal{B}\big(X(t)\big) \mbox{ for a.e.} \  t \in \left(0,\infty \right).
\end{equation}

To each solution $X(\cdot)$ of the differential inclusion, we can associate an extended control
$a(\cdot)=(\alpha_1,\alpha_2,\mu)(\cdot)\in A:=A_1\times A_2\times[0,1]$ so that $(X,a)$ becomes a
controlled trajectory of the system. 
%If the trajectories satisfy \eqref{t-bis}, we denote them by $X_x(\cdot)$. 
More precisely, introducing the following notations
$$\mathcal{E}_1:=\{t: \  X(t) \in \Omega_1  \}, ~\mathcal{E}_2:=\{t: \  X(t) \in \Omega_2 \},~ 
\mathcal{E}_{\mathcal{H}}:=\{t: \  X(t) \in \mathcal{H} \}\,,$$
then the result is the following. For more details see \cite[Th. 2.1]{manoel}.

% theorem
\begin{theorem}\label{dynamic}
     For each solution $X(\cdot)$ of the differential inclusion  \eqref{t}, 
     there exists a control $a(\cdot) =
    \left(\alpha_{1}(\cdot), \alpha_{2}(\cdot), \mu(\cdot)\right) \in \mathcal{A}:=  L^{\infty}_{loc}(0,\infty;A)$
    such that
    \begin{equation}\label{trajec} 
   \begin{split}
       \dot X(t) = &  ~b_1\big(X(t),\alpha_{1}(t)\big) \mathds{1}_{ \lbrace X(t)\in 
    \Omega_1\rbrace} + b_2\big(X(t),\alpha_{2}(t)\big) \mathds{1}_{ \lbrace X(t)\in 
    \Omega_2\rbrace} 
    \\
    & + b_{\mathcal{H}}\big(X(t),a(t)\big) \mathds{1}_{ \lbrace X(t)\in 
   \mathcal{H}\rbrace} \ \ \ \ \mbox{for  a.e. } \ \ t>0.
    \end{split}
    \end{equation}
    %%%%%%%%
    Moreover, 
    \begin{equation}\label{stampacchia} 
        b_{\mathcal{H}}\big(X(t),a(t)\big) \cdot e_N=0   \ \ \mbox{for a.e. } t\in\mathcal{E}_{\mathcal{H}}.
    \end{equation}
\end{theorem}

The set of trajectories starting from $x\in\R^N$ can then be defined by
\begin{equation}\label{control} 
    \tau_{A}(x):=\Big\{\big(X(\cdot),a(\cdot)\big) \in \mathrm{Lip}(\mathbb{R}^{+};
    \mathbb{R}^{n})\times \mathcal{A} \ \  \mbox{satisfying } \eqref{trajec}\ \ \mbox{and}\ X(0)=x
    \Big\}\,,
\end{equation}
and for any $x\in\H$, we define $A_0(x):=\big\{a\in\mathcal{A}: b_\H(x,a)\cdot e_N=0\big\}$,
the set of controls maintaining the trajectory on $\H$.

As we saw in the introduction, regular dynamics and trajectories are also of interest. They are
defined as follows:
\begin{definition}\label{def:regular}\
    \begin{itemize}
        \item[$(i)$] Given $x\in\H$ and $a=(\alpha_1,\alpha_2,\mu)\in A_0(x)$, we say that
            $a=(\alpha_1,\alpha_2,\mu)$ is a {\bf regular control} if 
            \begin{equation}\label{ineq:reg}
                b_1(x,\alpha_1)\cdot e_N \leq 0~~ \mbox{and} ~~b_2(x,\alpha_2)\cdot e_N \geq 0\,.
            \end{equation}
            We denote by $\Areg(x)$ the subset of controls in $A_0(x)$ which are regular.
        \item[$(ii)$] If $x\in\H$ and $\alpha\in\Areg(x)$, we say that the dynamic $b_\H(x,a)$
            is regular.
        \item[$(iii$)] The set of {\bf regular trajectories} is defined as
            \begin{equation}\label{controlregular} 
            \taureg(x):=\Big\{\big(X(\cdot),a(\cdot)\big) \in \tau_{A}(x): \mbox{ for a.e.} 
            \ t\in \mathcal{E}_{\mathcal{H}}\,, \ \ b_\H\big(X(t),a(t)\big)\mbox{ is regular}\Big\}\,.
            \end{equation}
        
    \end{itemize}
\end{definition}

We will also use some intermediate set between the set of regular and general trajectories in
Section~\ref{sect:U.eta}, by relaxing condition \eqref{ineq:reg} with a parameter $\eta>0$ instead
of $0$.

Notice that, given a controled trajectory $(X,a)$, we can define $(b,l)$
globally by setting $(b,l)=(b_i,l_i)$ in $\Omega_i$ for $i=1,2$ and $(b,l)=(b_\H,l_\H)$ on $\H$.
We then define two value functions by
\begin{equation}\label{umenos}
    U^{-}_{A}(x) := \inf_{(X,a)\in\tau_{A}(x)}\left( \int_{0}^{\infty} 
    l( X(t),a(t)) e^{-\lambda t} dt \right).
\end{equation}
\begin{equation}\label{umais}
    U_{A}^{+}(x) := \inf_{(X,a)\in\taureg(x)}\left( \int_{0}^{\infty} 
    l( X(t),a(t)) e^{-\lambda t} dt \right).
\end{equation}
We detail some properties of $U^-_A$ and $U^+_A$ in Section~\ref{sect:VF}. Of course, similar
definitions are used for the case of compact control sets, and the notation of the respective value functions 
are $U_{A^m}^-$ and $U_{A^m}^+$.  

%------------------------------------------------------------------------
\subsection{Hamiltonians}

As is well-known, optimal control problems are related to some Hamilton-Jacobi equation satisfied by
value functions. The natural Hamiltonian associated to the control problem is defined by 
\begin{equation}\label{h}
	H(x,u,p):=\sup_{a\in\mathcal{A}}\big\{\lambda u-b(x,a)\cdot p-l(x,a)\big\}\,,
\end{equation}
where we recall that the extended control takes the form $a=(\alpha_1,\alpha_2,\mu)$.
We refer to \cite[Lemma 7.3]{guy} for proof of the fact that
\begin{equation}\label{H.global} 
    H\left(x, u, p \right):= \left\{\begin{matrix} 
    H_1\left(x, u, p \right)  & \mbox{in }   \Omega_1,
    \\H_2\left( x, u, p \right)  & \mbox{in } \Omega_2, 
    \\ \max\{
    H_1\left(x, u, p \right),~H_2\left(x, u, p \right)\}  & \mbox{on }  \mathcal{H}.
     \end{matrix} \right.
\end{equation}

Moreover, as mentioned in the introduction, in order to get a comparison result for \eqref{ISHII}, we need to
define the tangential Hamiltonian. For $(x,u,q)\in \mathcal{H}\times\R\times \R^{N-1}$ we set
\begin{equation}\label{def:HT}
    H_{T}(x,u, q)=\sup\limits_{a\in A_0(x)} \Big\{ \lambda u
    - b_{\mathcal{H}}(x, a)\cdot (q,0) - l_{\mathcal{H}}( x,a)\Big\}\,.
\end{equation}
We can similarly define $H_T^\mathrm{reg}$ by using regular controls, as is done in
\cite{manoel}:
\begin{equation}\label{def:HTreg}
    H_{T}^\mathrm{reg}(x,u, q)=\sup\limits_{a\in A_0^\mathrm{reg}(x)}
    \Big\{ \lambda u - b_{\mathcal{H}}(x, a)\cdot (q,0) - l_{\mathcal{H}}( x,a)\Big\}\,,
\end{equation}
however, this second tangential Hamiltonian is not as useful as $H_T$ in the sense that it does not
lead to a characterization of $U^+$, nor a satisfying comparison result.
We now turn to the regularity of Hamiltonians, proving that the $H_i$ are well defined, continuous,
attain their supremum in a compact control set and are locally coercive with respect to the third
variable.

\begin{proposition}\label{H_i_cont}\label{reducao} 
 Assume that \HA, \HB and \HC hold. Then,
 \begin{enumerate}
     \item[$(i)$] If $K\subset \overline{\Omega_i}\times\R\times\R^N$ is compact, there exists
         $\tilde A\subset A_i$ compact such that for any $(x,u,p)\in K$, the supremum defining
         $H_i(x,u,p)$ is attained for $\alpha_i\in \tilde A$. 
     \item[$(ii)$] $H_i(x, u, p): \overline{\Omega_i}\times \mathbb{R} \times \mathbb{R}^{N}\to\R$
         is continuous with respect to all the variables. 
     \item[$(iii)$] $H_i(x,u,p)$ is coercive in $p$, locally uniformly with respect to $(x,u)$. 
\end{enumerate}
As a consequence, the same properties are valid for $H$ and $H_T$.
\end{proposition}

\begin{proof} 
    Recall that the Hamiltonian $H_i$ is defined by 
    $$H_i(x,u,p):=\sup_{\alpha\in A_i}\big\{\lambda u-b_i(x,\alpha)\cdot p-l_i(x,\alpha)\big\}\,.$$ 
 
    \noindent ({\it i}\,)
    On $X=\overline{\Omega_i}\times \mathbb{R} \times \mathbb{R}^{N}\times A_i$, let us 
    define $h_i: X \longrightarrow \mathbb{R}$ by
    $$h_i(x, u, p, \alpha_i):=\lambda u-b_i(x,\alpha_i)\cdot p -l_i(x,\alpha_i)\,,$$   
    and notice that $h_i$ is a continuous function. Notice also that there exists $M=M(K)>0$ such
    that $K\subset V_i\times[-M,M]\times B_M$, where $V_i:=\overline{\Omega}_i\cap B_M$, and $B_M$ 
    is a ball centered in zero with radius $M$. 
    Thanks to hypothesis \HB, for any
    $(x,u,p)\in K$, there exists a control $\alpha^x_i\in A_i^M\subset A_i$ such that
    $b_i(x,\alpha^x_i)=0$. 
    So, for any $(x,u,p)\in K$  
    $$\sup_{\alpha_i \in A_i}\big\{ -b_i(x,\alpha_i)\cdot p -l_i(x,\alpha_i)\big\} 
    \geq -b_i(x,\alpha_i^x)\cdot p -l_i(x,\alpha_i^x) \geq -\|l_i\|_{L^\infty(V_i\times A_i^M)}.$$
    Let us fix 
    \begin{equation}\label{Mi}
        \tilde M_i > \max\big\{\max_{|p|\leq M}|p|, \|l_i\|_{L^\infty(V_i\times A_i^M)}\big\}=
        \max\big\{M,\|l_i\|_{L^\infty(V_i\times A_i^M)}\big\}\,.
    \end{equation}
    Thanks to hypothesis \HC, there exists $\Gamma_{\tilde M_i} >0$ such that if $|\alpha_i| >
    \Gamma_{\tilde M_i} $ then
    \begin{equation}\label{bound_li}
	    l_i(x, \alpha_i) > \tilde M_i(1 + |b_i(x,\alpha_i)|) \quad  \forall x \in V_i.
    \end{equation}
    Hence, $- b_i(x,\alpha_i)\cdot p -l_i(x,\alpha_i)  <  -\|l_i\|_{L^\infty(V_i\times
    A_i^M)}$, for controls such that $|\alpha_i| \geq \Gamma_{\tilde M_i}$. Therefore, for any
    $(x,u,p)\in K$, the supremum of
    $h_i(x, u, p, \alpha_i) = \lambda u- b_i(x,\alpha_i)\cdot p-l_i(x,\alpha_i)$ is attained
    for $\alpha_i\in A_i^M$, which is a compact subset of $A_i$. This proves $(i)$.

    \smallskip

    \noindent ({\it ii}\,) The fact that $H_i$ is continuous just derives from $(i)$: since $h_i$ is
    continuous and the supremum is locally attained on a compact set of controls, the supremum of
    $h_i$ is continuous with respect to~$(x,u,p)$, which proves $(ii)$.

    \noindent ({\it iii}\,) Let us turn to the coercivity property.
    Let  $(x,u)\in V_i\times[-M,M]$ and $p \in \mathbb{R}^{N}$. Due to \HB, for any $x\in V_i$,
    there exists $\tilde\alpha_i^x\in A_i^M $ such that $|b_i(x,\tilde\alpha_i^x)|=\delta$ and
    $-b_i(x,\tilde\alpha_x)\cdot p =\delta|p|$. This implies that
    $$H_i(x,u,p)  \geq \lambda u +\delta\vert p \vert -||l_i||_{L^\infty(V_i\times A_i^M)}\,,$$ 
    therefore, $\lim_{|p|\rightarrow \infty }H_i(x,u,p)=+\infty$, and $(iii)$ holds.

    \smallskip

    Looking at the definition of $H$ and $H_T$, it is clear
    that those Hamiltonians enjoy the same properties as $H_1$ and $H_2$, the proofs being
    essentially identical.
\end{proof}

\

\subsection*{Approximations}

Let us introduce a family $(H^m)_{m\in\N_*}$ of Hamiltonians by reducing the control problem to
compact control sets $A^m=A_1^m\times A_2^m\times[0,1]$: 
\begin{equation}\label{Hm}
	H^m(x,u,p):=\sup_{a\in\mathcal{A}^m}\big\{\lambda u-b(x,a)\cdot p-l(x,a)\big\}\,,
\end{equation}
where of course $\mathcal{A}^m:=L^\infty_{loc}(0,\infty;A^m)$. It is clear from the definition of $A^m$ that
for any fixed $(x,u,p)$, $H^m(x,u,p)$ converges monotonically to $H(x,u,p)$.
Notice that both Hamiltonians are upper semi-continuous by construction, and that of course for
fixed $(x,u,p)$, $H^m(x,u,p)\to H(x,u,p)$ monotonically.

Now, in order to connect $H^m$ and $H$ we need to introduce the half-relaxed limits of $H^m$.
We refer the reader to \cite{guy} for more precise results on the semi-continuous enveloppes as well
as upper and lower semi-continuous functions (usc and lsc for short). Let us just recall here that for any
locally bounded function $f$ on a set $\Omega\subset\R^N$, the lower and upper semi-continuous
enveloppes are defined respectively by 
$$f_*(x):=\liminf_{y\to x}f(y)\,,\quad f^*(x):=\limsup_{y\to x}f(y)\,.$$
Moreover, if $(f_\e)_{\e>0}$ is a family of locally uniformly bounded functions, the half-relaxed
limits as $\e\to0$ are defined similarly as
$$\liminfs_{\e\to0} f_\e(x):=\liminf_{\e\to0\atop y\to x}f_\e(y)\,,\quad
\limsups_{\e\to0} f_\e(x):=\limsup_{\e\to0\atop y\to x}f_\e(y)\,.$$
Of course, $\liminfs f_\e$ and $\limsups f_\e$ are respectively lower semi-continuous and upper
semi-continuous. 
We will need the following result, which appears in \cite{guy} as an exercise (we provide here a
full proof for completeness):
\begin{lemma}\label{lem:half.relaxed}
    Let $(u_\e)_{\e>0}$ be a non-decreasing (with respect to $\e$)
    sequence of continuous functions which are locally uniformly
    bounded on some set $\Omega\subset\R^N$. Then,
    $$\begin{aligned}
        \liminfs u_\e &= \sup_{\e>0} u_\e\,,\\
        \limsups u_\e &= \big(\sup_{\e>0} u_\e\big)^*\,.
    \end{aligned}$$
\end{lemma}
\begin{proof}
    Let us first notice that since $(u_\e)$ is locally uniformly bounded, the various quantities in
    this result are well-defined. Notice also that by monotonicity, 
    as $\e\to0$, $u_{\e}(x)\to \big(\sup_{\e>0} u_\e\big)(x)$.
    
    \smallskip

   Let $\e'>0$ be fixed and consider a sequence $(\e_n,y_n)\to (0,x)$. Then for
    $\e_n<\e'$, \ie, $n$ big enough,
    $$\begin{aligned}
        u_{\e_n}(y_n) &= \Big[u_{\e_n}(y_n)-u_{\e'}(y_n)\Big]+
        \Big[u_{\e'}(y_n)-u_{\e'}(x)\Big]+u_{\e'}(x)\\
        &\geq o_n(1)+u_{\e'}(x)\,.
        \end{aligned}$$
    Indeed, this is obtained by using the monotonicity of $(u_\e)$ for the first bracket, 
    and the continuity of $u_{\e'}$ for the second one.
    By taking the liminf, we deduce that $\liminfs u_\e(x)\geq u_{\e'}(x)$ for any $\e'>0$. Now,
    since by monotonicity, $u_{\e'}\to \sup u_\e$, we get the inequality
    $\liminfs u_\e \geq \sup_{\e>0} u_\e\,.$
    For the converse inequality, it is enough to notice that taking $(\e_n,x)\to(0,x)$ leads to
    $$\liminfs u_\e(x)\leq \liminf_{\e_n\to0}u_{\e_n}(x)=\sup_{\e>0} u_\e(x)\,$$
    and the first result of the Lemma follows.

    \smallskip

    For the limsup property, we first notice that of course
    $\limsups u_{\e} \geq \liminfs u_{\e}=\sup_{\e>0} u_{\e}\,,$
    and since $\limsups u_\e$ is upper semi-continuous, necessarily
    $\limsups u_\e\geq \big(\sup_{\e>0} u_\e\big)^*\,.$
    For the converse inequality, using the monotonicity property of $(u_\e)$ we see that
    $$\limsups u_\e (x)\leq \limsup_{y\to x}\, \big(\sup_{\e>0} u_\e\big)(y)=\Big(\sup_{\e>0}u_\e\Big)^*$$
    and the second result follows.
\end{proof}

Let us give and important consequence:
\begin{corollary}\label{cor:Hm.UAm}
    The following limits hold:
    $$\limsups_{m\to+\infty} H^m(x,u,p) = H(x,u,p)\,.$$
    $$\liminfs_{m\to+\infty} U_{A^m}^-(x)=U_A^-(x)\,.$$
\end{corollary}
\begin{proof}
    This is a direct consequence of the previous Lemma. Take $u_{\e_n}=H^m$: we already noticed that
    for fixed $(x,u,p)$, the sequence $(H^m(x,u,p))_m$ is non-decreasing; each $H^m$ is continuous
    and $H^m\leq H$ which yields a local uniform bound.  Hence, using that $H$ is upper
    semi-continuous,
    $$\limsups_{m\to+\infty}H^m(x,u,p)=\Big(\sup_m H^m(x,u,p)\Big)^*=
    \Big(H(x,u,p)\Big)^*=H(x,u,p)\,.$$
    For the case of $U_{A^m}^-$, let us notice that since $A^m$ is compact, $U_{A^m}^-$ is
    continuous. Moreover, $m\mapsto U_{A^m}^-$ is non-increasing, and that
    $$\begin{aligned}
        \lim_{m\to+\infty}U_{A^m}^-(x) &= \inf_{m\in\N_*} U^{-}_{A^m}(x)= 
        \inf_{m\in\N_*}\inf_{(X,a)\in\tau_{A^m}(x)}\left( \int_{0}^{\infty} l( X(t),a(t)) e^{-\lambda
        t} dt \right)\\
        &= \inf_{(X,a)\in\tau_{A}(x)}\left( \int_{0}^{\infty} l( X(t),a(t)) e^{-\lambda
        t} dt \right)= U_A^-(x)\,.
    \end{aligned}$$
    Then, by using again Lemma~\ref{lem:half.relaxed} with the nondecreasing sequence
    $(-U_{A^m}^-)$, we end up with the second result.
\end{proof}

%------------------------------------------------
\subsection{Viscosity solutions}

By solutions of Problem~\ref{ISHII} we mean of course \textbf{viscosity solutions}. Let us briefly
recall some definitions:
\begin{definition}\label{def:viscosity}
    Let us consider an abstract Hamilton-Jacobi equation $H(x,u,Du)=0$ posed in a set $\Omega$.
    \begin{itemize}
        \item[$(i)$] An upper semi-continuous function $u:\Omega\to\R$ is a viscosity subsolution in
            $\Omega$ if for any $C^1$-smooth test-function $\phi$, at any local maximum point
            $x\in\Omega$ of $u-\phi$ we have
            $$H(x,u,D\phi(x))\leq0\,.$$
        \item[$(ii)$] A lower semi-continuous function $v:\Omega\to\R$ is a viscosity supersolution
            in $\Omega$ if for any $C^1$-smooth test-function $\phi$, at any local minimum point
            $x\in\Omega$ of $u-\phi$ we have
            $$H(x,u,D\phi(x))\geq0\,.$$
        \item[$(iii)$] A continuous function $u:\Omega\to\R$ is a viscosity solution of
            $H(x,u,Du)=0$ in $\Omega$ if it is at the same time a subsolution and a supersolution.
    \end{itemize}
\end{definition}

In the rest of the paper, we use the abbreviation \usc and \lsc for upper semi-continuous and lower
semi-continuous respectively.

Notice that if $u$ is only (locally) bounded, we can still consider its upper semi-continuous and
lower semi-continuous enveloppes, $u^*$ and $u_*$, which allows to define sub/super solutions through $u^*$
and $u_*$. But in this paper, we will always consider \usc subsolutions and \lsc supersolutions.

The Ishii conditions in \eqref{ISHII}, \ie the min and max conditions on $\H$ are checked as above
by using also $C^1$-smooth test functions $\phi$, where $(i)$ $u-\phi$ reaches a maximum at $x\in\H$ for
the min equation (subsolution condition); $(ii)$ $u-\phi$ reaches a minimum at $x\in\H$ for the max
equation (supersolution condition).

The case of $H_T$ or $H_T^\mathrm{reg}$ is particular since it is a $(N-1)$-dimensional equation. So,
Definition~\ref{def:viscosity} has to be applied with test-functions $\phi:\H\to\R$ and considering
local maximum/minimum of $u-\phi$ on $\H$. Equivalently, we can use $\phi$ as a test function in
$\R^N$ and consider the max/min in $x'\in\H$ of $u(x',0)-\phi(x',0)$.

We end up this section by defining the notion of stratified solution, by adding $H_T$ to the set of
inequations. Following \cite{Barles2018} we introduce the following definition.
\begin{definition}\label{def:stratified}
    Let us consider problem \eqref{ISHII} that we write under the abstract form $\HH(x,u,Du)=0$.
    \begin{itemize}
        \item[$(i)$] An \usc function $u:\R^N\to\R$ is a stratified subsolution of 
            \,$\HH(x,u,Du)\leq0$ $\HH(x,u,Du)=0$ if it is a subsolution of \eqref{ISHII}
            satisfying the additional inequality $H_T(x,u,Du)\leq 0$ in $\H$.
        \item[$(ii)$] A \lsc function $v:\R^N\to\R$ is a stratified supersolution of 
            \,$\HH(x,u,Du)\geq0$ $\HH(x,u,Du)=0$ if it is a supersolution of \eqref{ISHII}.
        \item[$(iii)$] A continuous function is a stratified solution of \,$\HH(x,u,Du)=0$ if it is at
            the same time a stratified subsolution and a stratified supersolution.
    \end{itemize}
\end{definition}
Notice that of course, only the subsolution condition is complemented by $H_T\leq0$, stratified
supersolutions are nothing but standard (Ishii) supersolutions of \eqref{ISHII}.
We also want to comment on the fact that in \cite{Barles2018}, the notion of weak and strong
stratified solutions are used, due to some specific regularity properties of subsolutions. In our
situation, we are assuming total controllability and Ishii inequalities on $\H$ so that weak
stratified subsolutions are necessarily also strong ones. For more details, see \cite[Prop. 19.2]{Barles2018}.

One of the main goals of this paper is to establish a global comparison result for stratified
solutions of $\HH(x,u,Du)=0$, or more precisely problem \eqref{ISHII}. By global and local
comparison results, as in \cite{Barles2018} we mean here:

\

\noindent(GCR): \textit{For any \usc subsolution, $u$, and \lsc supersolution, $v$, $u\leq v$ in $\R^N$.} 

\

\noindent(LCR): \textit{For any $x\in\R^N$, there exists $\bar r>0$ such that for any 
\usc subsolution, $u$, and \lsc supersolution, $v$, and any $0<r<\bar r$,}
$$\max_{\overline{B_r(x)}}(u-v)^+=\max_{\partial B_r(x)}(u-v)^+\,.$$

\

We refer to Section~\ref{sect:comp} below where both type of results are proved.

%%%%%%%%%%%%%%%%%%%%%%%%%%%%%%%%%%%%%%%%%%%%%%%%%%%%%
\section{Value functions}\label{sect:VF}
%%%%%%%%%%%%%%%%%%%%%%%%%%%%%%%%%%%%%%%%%%%%%%%%%%%%%

Let us begin by recalling that in the bounded control case, \ie, considering control sets $A^m$ for
$m\in\N_*$, the results in \cite{manoel} directly apply. If $H_i^m$, $H_T^m$ and
$H_T^{m,\mathrm{reg}}$ are the Hamiltonians associated to the compact control set $A^m$ and
$U^-_{A^m},\ U^+_{A^m}$ are the associated value functions, the following result holds:

\begin{theorem}[\cite{manoel}, Thm. 2.5] \label{existenciam}
    Assuming \HA, \HB and \HC, $U^{+}_{A^m}$ and $U^{-}_{A^m}$ are viscosity solutions of the Ishii
    problem \eqref{ISHII}. 
    \begin{enumerate}
        \item[$(i)$] $U^{-}_{A^m}$ is a stratified subsolution, associated to the tangential
            Hamiltonian $H^m_T$, \ie, for any $\bar{x}\in\R^{N-1}$, $U^{-}_{A^m}(\bar{x},0)$
            satisfies 
            $$ H^m_T(x,u, D_{\mathcal{H}}u)\le 0\,.$$
        \item[$(ii)$] $U^{+}_{A^m}$ is a supersolution of the tangential Hamiltonian
            $H^{m,\mathrm{reg}}_T$, \ie, for any $\bar{x}\in\R^{N-1}$, $U^{+}_{A^m}(\bar{x},0)$
            satisfies $$ H^{m,\mathrm{reg}}_T(x,u, D_{\mathcal{H}}u)\ge 0\,.$$ 
    \end{enumerate}
\end{theorem} 

We now go back to the non-bounded case.

%--------------------------------------------------------------------
\subsection{The dynamic programming principle}

The result below is the Dynamic 
Programming Principle which is the main result to prove that the value functions are subsolutions
and supersolutions of \eqref{ISHII}.
% theorem
\begin{theorem}[Dynamic Programming Principle] \label{dppt}
\begin{equation}\label{prdinmenos}
    U_{A}^{-}(x)  = \inf_{(X,a)\in\tau_{A}(x)}\left( \int_{0}^{T} l( X(t),a(t)) e^{-\lambda t} dt +
    U_{A}^{-} \left(X\left(T \right)  \right)e^{-\lambda T}\right).  
\end{equation}

\begin{equation}\label{prdinmais}
    U_{A}^{+}(x)  = \inf_{(X,a)\in\tau_{A}^{\mathrm{reg}}(x)}\left( \int_{0}^{T} l( X(t),a(t))
    e^{-\lambda t} dt + U_{A}^{+}\left(X\left(T \right)  \right)e^{-\lambda T}\right).  
\end{equation}
\end{theorem}

\begin{proof}
    The proof is standard, cf.  \cite[p. 65]{guy}.
\end{proof}

In the following result, we prove that the trajectories in $\tau_{A^m}(x)$ are locally bounded for time small enough. 
\begin{lemma}\label{trjbounded}
Let $m\in\mathbb{N}*$ and let $b_i:\overline{\Omega_i}\times A^m_i\to\mathbb{R}^N$ be continuous. Given 
$r>0$ and $x \in \mathbb{R}^{N}$, there exist $\overline{t}>0$  such that for any trajectory $(X,a) \in 
\tau_{A^m}(x)$, we obtain
$$|X_x(s)-x|\leq r,  ~\forall~ s\leq \overline{t} . $$
\end{lemma} 

\begin{proof}
    Integrating \eqref{trajec}, we consider $\overline{t} \leq r\cdot(||b||_{L^\infty(B_r(x)\times
    A^{m}})^{-1}$ to obtain the result.
\end{proof}

%--------------------------------------------------
\subsection{Regularity}
  
The main result here is the following
\begin{proposition}\label{continuouss}
    Let us assume that hypotheses \HA, \HB, \HC hold.  Then $U_A^{-} $ and $
    U_A^{+}$ are locally bounded and locally Lipschitz continuous functions.  
\end{proposition}

\begin{proof} We first prove some local bound, the Lipschitz regularity follows almost directly.

    \smallskip

    Let $V \subset \R^N$ be a convex compact set, and assume $x\in \Omega_i\cap V$.  Thanks to
    hypothesis \HB, there exists a control $\alpha_i^*\in A_i$ such that $b_i(x,\alpha^*_i)=0$,
    yielding the constant trajectory $X_x(t)=x$ with associated control $\alpha_{i}(t)=
    \alpha_i^{*}$ for all $t$. Hence, 
    \[
        \begin{array}{ll} 
        U^-_A(x) &\displaystyle =\inf_{\tau(x)}\left( \int_{0}^{\infty} l(
        X_x(t),\alpha(t)) e^{-\lambda t} dt \right) \le \int_{0}^{\infty} l( x,\alpha_i^*)
            e^{-\lambda t}dt=\frac{l( x,\alpha_i^*)}{\lambda} \,.
        \end{array}
    \]
    Using hypothesis \HA, it is clear that $U^-_A$ is bounded in $\Omega_i \cap V$.

    Now, if $x \in \mathcal{H}\cap V$, we use the same approach with some controls $\alpha_1^{*} \in
    A_1$, $\alpha_2^{*} \in A_2$ such  that  $b_1(x,\alpha_1^{*})= 0 = b_2(x,\alpha_2^{*})$. The
    trajectory $X(t)=x$ associated to $a(t)= (\alpha_1^{*}, \alpha_2^{*}, \mu)$, with
    $\mu=1/2$ for instance, yields a similar bound. So, $U^-_A$ is also bounded in $V\cap\H$ and
    finally, $U^-_A$ is bounded in $V$. Of course the same arguments work for $U^+_A$ because the
    control $a=(\alpha_1^*,\alpha_2^*,\mu)$ is obviously regular.
  
    \smallskip

     For the local Lipschitz continuity, the proof is analogous to \cite[Th.
    2.3]{manoel}: fix as above $V\subset\R^N$ compact and given $x,y\in V$, let us consider the
    constant dynamic
    $$\tilde b:=-\delta \frac{x-y}{|x-y|}\,.$$
    Since $|\tilde b|=\delta$, by \HB,  $\tilde b\in\mathcal{B}^{A_i^m}(z)$ for all $z\in[x,y]$, where
    $m\in\N_*$ is fixed. Notice here
    that $x,y$ may be on the same side of $\H$ or on the opposite sides, and even one or both may be
    located on $\H$. In all these cases, the trajectory $X(t):=x+t\tilde b $ is an admissible
    straight line associated to some extended control $a(\cdot)\in A^m$, \ie $|a(\cdot)|\leq m$, such that
    $$X(0)=x\,,\quad X(|x-y|/\delta)=y\,.$$
    Using $X(t)$ in the dynamic programming principle yields
    $$U^-_A(x)\leq \int_0^{{|x-y|\over\delta}}l(X(s),a(s))e^{-\lambda s}ds + U^-_A(y)e^{-\lambda {|x-y|\over \delta}}\,.$$
    Using \HA, we deduce that
    $$U^-_A(x)-U^-_A(y)\leq \Big(\sup_{V\times A^m}l(z,a)\Big){|x-y|\over \delta}+\Big(e^{-\lambda {|x-y|\over \delta}}-1\Big)
    \sup_V U^-_A(z)\leq  C|x-y|$$
    for some constant $C=C(V)$. Since $x,y$ are arbitrary in $V$ this implies the local Lipschitz
    bound: $|U_A^-(x)-U_A^-(y)|\leq C(V)|x-y|$. 

    The same arguments work for $U_A^+$ because the trajectory is of course regular: if $x,y$ are
    located on the same side of $\H$ there is nothing to do; if they are on opposite sides the
    trajectory is just reaching $\H$ for one specific time $t$ but we do not take it into account
    since this is a neglectable set in time; finally if $x,y\in\H$ the trajectory is purely
    tangential by construction, so it is regular.
  \end{proof}

%-------------------------------------------------------------------
\subsection{Value functions are solutions of the Ishii problem}

In order to prove that value functions $U_A^-$ and $U_A^+$ are both solutions in the sense of Ishii of
problem~\eqref{ISHII}, we use the approximations by compact control sets and pass to the limit.

\begin{theorem}\label{existencia}
    Assume that \HA, \HB and \HC hold. Then 
    \begin{enumerate}
        \item[$(i)$] the value functions $U_A^{-}$ and $U_A^{+}$ are both viscosity solutions
            of~\eqref{ISHII} in the sense of Ishii;
        \item[$(ii)$] $x'\mapsto U_A^{-}(x',0)$ satisifies the tangential subsolution inequality 
            $$H_T(x,u,Du)\leq0\quad\text{on}\quad \H\,.$$
    \end{enumerate}
 \end{theorem}

\begin{proof} We first prove the supersolution property, then the subsolution one and finally the
    tangential property of $U^-_A$.

    \smallskip

    \noindent ({\it i}\,) 
    In order to prove that $U_A^{-}$ is a supersolution, we proceed as in
    \cite{unbounded}. We recall that by
    Corollary~\ref{cor:Hm.UAm}, the following limits hold:
    $$\limsups_{m\to+\infty} H^m(x,u,p) = H(x,u,p)\,,\quad
    \liminfs_{m\to+\infty} U_{A^m}^-(x)=U_A^-(x)\,,$$
    which allow to apply the half-relaxed limit method directly (cf. \cite[Section 2.1.2]{Barles2018}): 
    since for each $m$ fixed,
    $U_{A^m}^-$ is a supersolution of the Ishii problem associated to $H^m$, then 
    $U_A^-=\liminfs U_{A^m}^-$ is also a supersolution of the problem associated to $H=\limsups
    H^m$. In view of \eqref{H.global}, in other words, we have proved that $U_A^-$ is a
    supersolution of \eqref{ISHII}. Of course the same argument is valid for $U_A^+$.

    \smallskip

    \noindent ({\it ii}\,) To prove that are $U_A^{+}$ and  $U_A^{-}$ are subsolutions we
    proceed exactly as in \cite[Thm. 2.5]{manoel}. For the sake of completeness we provide here the
    main arguments but the proof readily applies to our case, using some compactness arguments that
    we proved above (typically, the fact that the $H_i$ are continuous, etc.).
    Of course the proof in each $\Omega_i$ is standard so
    we focus on getting the inequality $\min(H_1,H_2)\leq0$ on $\H$, and we only do so for $U_A^+$
    since the proof for $U_A^-$ is similar.

    Let $\phi \in C^{1}(\R^N)$ and consider $x \in \mathcal{H}$, a local maximum point of
    $U^+_A-\phi$. Assuming without loss of generality that this maximum is zero, there exists $r>0$
    such that $U_A^+(y)-\phi(y)\leq 0 $ for all  $y \in B_r(x) \subset \R^{N}$ and
    $U_A^+(x)=\phi(x)$. Thanks to Lemma \ref{H_i_cont}, $H_i$ attains its
    supremum in a bounded control set, so there exist $(\alpha_1,\alpha_2)\in A_1\times A_2$ such
    that 
    \[
    	\begin{array}{l}
    	H_1\left(x, \phi\left( x\right), \nabla\phi \left( x\right) \right)= \lambda \phi\left( x \right) 
    	-b_1\left(x,\alpha_1 \right)\cdot\nabla\phi \left( x\right)-l_1(x,\alpha_1).
    	\\
    	H_2\left(x, \phi\left( x\right), \nabla\phi \left( x\right) \right)  =  \lambda \phi\left( x \right) 
    	-b_2\left(x,\alpha_2 \right)\cdot\nabla\phi \left( x\right)-l_2(x,\alpha_2).
    	\end{array}
    \]
    In order to prove that $U_A^+$ is subsolution of $\min\{ H_1,H_2\}$, we use specific
    trajectories that we build using the constant controls $\alpha_1$ and $\alpha_2$. Notice that
    such controls may not necessarily be regular, but in each case we find a suitable regular
    trajectory in order to use the dynamic programming principle for $U_A^+$. 

    Let us focus on the (regular) case where 
    $$ b_1\left(x,\alpha_1 \right) \cdot e_N < 0 \,\mbox{ and }\,
    b_2\left(x,\alpha_2 \right) \cdot e_N > 0.$$ 
    The proof for all the other posible combinations of signs of $b_i(x,\alpha_i)\cdot e_N$ is
    similar with a few adaptations using normal controllability, see \cite[Thm. 2.5]{manoel}. 
    
    The idea is to construct a regular trajectory, $ (X,a)\in \tau^{\mathrm{reg}}_{A}(x)$, staying
    on $\mathcal{H}$, at least for a while.  Since the dynamic functions $b_i$ are continuous,
    there exists $\bar{\delta}>0$, such that for any $y\in\H\cap B_{\bar{\delta}}(x)$ the following
    quantity is well-defined: 
    \begin{equation}\label{mu_ca_A1} 
        \mu(y):= \frac{-b_2(y,\alpha_2)\cdot e_N}{(b_1(y,\alpha_1)-b_2(y,\alpha_2))\cdot
        e_N}\in(0,1)\,.
    \end{equation}
    Next, we consider the local trajectory defined by solving for $t>0$ small enough,
    \begin{equation}\label{trajA1}
    	\dot X(t) = \mu\big(X(t)\big) b_1\big(X(t),\alpha_1\big) 
    	+\big (1-\mu (X(t))\big)b_2\big(X(t), \alpha_2)\,,
    \end{equation}
    with $X(0)=x$.
    Since $b_1$, $b_2$ and $\mu$ are continuous functions, \eqref{trajA1} has a local solution and
    substituting $ \mu(X(t))$ in \eqref{trajA1}, we easily check that by construction, 
    $\dot X(t) \cdot e_N=0$ for $t>0$ small enough. Hence $X(\cdot)$  stays on $\mathcal{H} $ for a 
    while. Moreover, thanks to the continuity of $X(\cdot)$, there exists $T>0$ such that 
     $$
     	0 < \mu(X(t))< 1,\ \  b_1(X(t),\alpha_1)\cdot e_N<0, \ \ \ b_2(X(t),\alpha_2)\cdot e_N>0\ \ \ \mbox{for}  \ \ \ 
    	0\leq t\leq T.
    $$
    Let us define the following trajectory 
    \begin{equation}\label{traj0}
     	X_x(t):= \left\{
	\begin{matrix} 
    	 X(t), &   \mbox{if } 0 \leq t < T
    	\\  
	X(T), & \mbox{if }  t \geq T.
    	 \end{matrix} \right.
    \end{equation} 
    By the continuity of $b_i$ and the trajectory, the trajectory is Lipschitz.
    From hypothesis \HB, there exists $\alpha^{*}_i \in A_i $ such that  
    $b_i(X_x(T), \alpha_i^{*})=0$, and we may consider the extended control
     \begin{equation}\label{traj}
     a (t):= \left\{\begin{array}{ll} 
      (\alpha_1,\alpha_2,\mu(X(t)))&   \mbox{if } 0 \leq t < T,
    \\  (\alpha_1^{*}, \alpha_2^{*}, \mu) & \mbox{if }  t \geq T,
     \end{array} \right.
    \end{equation} 
     obtaining a regular trajectory $(X_x,a) \in \tau_{A}^{\mathrm{reg}}$. 
     Observe that $\mu$ for $t\geq T$ is any value in $[0,1]$.
    
    Recall $U_A(y)-\phi(y)\leq 0 $ for all $y \in B_r(x) $ and $U_A(x)=\phi(x)$. Since $X_x$ is
    Lipschitz, then there exists $0<T'<T$ such that $ X_x\left(  t \right) \in B_r(x)$ for all  
    $t < T'$. Now, by the Dynamic Programming Principle, we have
    $$
        \phi(x)\leq\int_{0}^{T'} l_{\mathcal{H}}\big(X_x(t),a(t)\big) e^{-\lambda t} dt + 
        \phi\big(X_{x}(T)\big)e^{-\lambda T'}.
    $$
    From the Fundamental Calculus Theorem, we have 
    $$ 
        0 \geq  \int_{0}^{T'}\Big[ -l_{\mathcal{H}}\big( X_x(t),a(t)\big) + \lambda\phi\big(X_x(t)\big) 
         -b_{\mathcal{H}}\big(X_{x}(t),a(t)\big)\cdot\nabla\phi \big(X_x(t)\big) \Big] e^{-\lambda t} dt.
    $$
    Dividing by $T'$ and taking the limit as $T'$ goes to $0$, we get
    $$ 
    \begin{aligned}
    0 & \geq \lambda \phi(x) - b_{\mathcal{H}}\big(x,\alpha_1, \alpha_2, \mu(x)\big)\cdot
        \nabla\phi(x)-l_{\mathcal{H}}\big( x,\alpha_1, \alpha_2, \mu(x)\big)\\ 
        &= \mu(x)\big(\lambda \phi(x)-b_1(x,\alpha_1)\cdot\nabla\phi(x)-l_1(x,\alpha_1)\big)\\ 
        &+\big(1-\mu(x)\big)\big(\lambda
        \phi(x)-b_2(x,\alpha_2)\cdot\nabla\phi(x)-l_2(x,\alpha_2)\big). 
    \end{aligned}$$
    So, finally we end up with,
    $$
    \begin{aligned}
    0 &\geq \min \big\{ \lambda \phi(x)-b_1(x,\alpha_1)\cdot\nabla\phi(x)-l_1(x,\alpha_1),
        \lambda \phi(x)-b_2(x,\alpha_2)\cdot\nabla\phi(x)-l_2(x,\alpha_2) \big\}  \\
        & =\min \big\{ H_1\big(x,\phi(x), \nabla\phi(x)\big),  H_2\big(x,\phi(x),
        \nabla\phi(x)\big)\big\}\,.
    \end{aligned}
    $$ 
    As we said, the other cases are treated by the same approach, the conclusion being that
    $U^+_A$ is subsolution of $\min \{H_1,H_2\}$ on $\mathcal{H}$. 
    
    \smallskip
    Proving the tangential inequalities is analogous: by constructing a
    trajectory staying on $\H$, the extended control $a(\cdot)$ belongs to $A_0$ and we get the
    subsolution inequality $H_T(x,u,Du)\leq0$ as above, by using the dynamic programming and passing
    to the limit. We refer again to \cite[Thm. 2.5]{manoel} for more details.
\end{proof} 

%%%%%%%%%%%%%%%%%%%%%%%%%%%%%%%%%%%%%%%%%%%%%%%%
\section{Comparison results}\label{sect:comp}
%%%%%%%%%%%%%%%%%%%%%%%%%%%%%%%%%%%%%%%%%%%%%%%%

In this section, we prove that stratified subsolutions and supersolutions of the Ishii problem are ordered.
We start giving a local comparison result for the Ishii problem, meaning that the maximum of $u-v$
in the closure of a ball is attained in the boundary. Then, under some localization hypotheses we
deduce a Global Comparison Result in $\R^N$. As a consequence of this comparison results, the value 
functions $U_A^-$ and
$U_A^+$ are extremal Ishii solutions. 

We recall that by definition, a stratified subsolution satisfies the additional subsolution
inequality $H_T\leq0$ while a stratified supersolution is nothing but a usual Ishii supersolution.

%-----------------------------------------------------------
\subsection{Local comparison result}\label{subsect:lcr}

Let us first recall that in the case of compact control sets, a Local Comparison Result for
stratified solutions can be found in \cite{manoel2014}, which translates here directly as a Local
Comparison Result for each problem associated to the control set $A^m$:
\begin{theorem}[see \cite{manoel2014}]\label{tcompracaoppp}
    Assume hypotheses \HA, \HB and let $m\in\N_*$ be fixed. Let $u$ be a locally bounded stratified subsolution
    of Ishii problem \eqref{ISHII} with Hamiltonians $H_1^m$,$H_2^m$,$H_T^m$ and $v$ be a locally
    bounded supersolution of the same problem. Then for any
    $x\in\R^N$ and $r>0$, the following comparison result holds:
    $$\max_{B_{r}(x)}(u-v)^{+} \leq \max_{\partial B_{r}(x)}(u-v)^{+}.$$
\end{theorem}

While it is clear from the definition that if $H_i^m(x,u,Du)\leq0$ in $\R^N$ the same holds for
$H_i\leq H^m_i$, getting a similar property for supersolutions can only be obtained under some
restrictions. In order to do so, let us introduce the notation $V\subset\subset\R^N$, meaning that
$V\subset\R^N$ is open, and $\overline{V}$ is compact. The result is the following.
\begin{proposition}\label{tcompracao}
    Assume hypotheses \HA, \HB, \HC and let $v$ be a locally bounded \lsc supersolution 
    of Ishii problem \eqref{ISHII}. Then, for any $V\subset\subset\R^N$,  
     there exists $m=m_V\in\N_*$ such that $v$ is a supersolution in $V$ 
    of the Ishii problem associated to $A^{m}$.
\end{proposition}

\begin{proof}
    First we fix $\bar{m}\in\N_*$ arbitrary. Due to hypothesis \HA, the following quantity is well-defined:
    $$K:=\frac{ ||l||_{L^\infty(V \times A^{\bar{m}}_i)} +\lambda ||v||_{L^\infty(V)} }{\delta}\,.$$
    To begin, we fix a test function $\phi \in C^1(V)$ such that $x\in V\cap\Omega_i$ is a local
    minimum point of $v-\phi$, (we explain how to treat the case $x\in\H$ at the end of the proof). 
   
    \smallskip

     Let us assume on one hand, that $|\nabla \phi (x)|\geq K$. 
    Notice that thanks to \HB, for any $\omega \in \overline{B_{\delta}(0)}$ and $x \in
    V\cap\Omega_i$, there exists
    $\alpha_i \in A^{\bar{m}}_i$ such that $b_i(x,\alpha_i)=\omega$. Applying this to 
    $\omega:=-\delta \nabla \phi (x)/\vert\nabla\phi(x)\vert$,  
    we find some $\alpha_i \in A^{\bar{m}}_i$ such that $b_i(x,\alpha_i)=\omega$.

    It follows that $-b_i(x,\alpha_i)\cdot\nabla\phi(x) = \delta\vert \nabla \phi (x)\vert$, which
    implies that 
    $$H^{\bar{m}}_i\left(x,v(x),\nabla \phi (x)\right) \geq -\lambda ||v||_{L^\infty(V)}
    +\delta\vert \nabla \phi (x) \vert - ||l_i||_{L^\infty(V \times A^{\overline{m}}_i)}  \geq 0.$$
  
    \smallskip

    Let us assume now that on the other hand, $|\nabla \phi (x)| < K$. Then,
    thanks to Proposition~\ref{reducao}, there exists $m(K)>0$ independent of $x\in V\cap \Omega_i$ such that
    $$H^{m(K)}_i(x,v(x),\nabla \phi (x))=H_i(x,v(x),\nabla \phi (x)) \geq 0.$$

    \smallskip

    Taking $m_V:=\max\{\bar m,m(K)\}$, we have that for any $x\in V\cap\Omega_i$, 
    $$H^{m_V}_i(x,v(x),\nabla \phi (x))\geq 0\,,$$
    because the supremum is taken over a larger set in each case. 
    Therefore, $v$ is a supersolution of $H^{m_V}_i$ in $\Omega_i$.

    \smallskip

    Let us study now the case, when $x \in \mathcal{H}$. 
    Of course, if $x \in V \cap \mathcal{H}$, then since
    $$\max \{H_1(x,v(x),\nabla \phi (x)), H_2(x,v(x),\nabla \phi (x)) \}\geq 0\,,$$
    there exists $i \in \{1, 2 \}$ with $H_i(x,v(x),\nabla \phi (x))\geq 0$ and from the above
    arguments we deduce that for some $m_V\in\N_*$,
    $$\max \{H^{m_V}_1(x,v(x),\nabla \phi (x)), H^{m_V}_2(x,v(x),\nabla \phi (x)) \}\geq 0\,,$$
    which ends the proof.
\end{proof}
 
Now we can prove the local comparison result for the Ishii problem associated to the complete
problem, associated to the unbounded control set $\mathcal{A}$.  
\begin{corollary}\label{tcompracaolllLL}
    Assume \HA, \HB and \HC. Let $u$ be a locally bounded stratified subsolution of Ishii problem
    \eqref{ISHII} and let $v$ be a locally bounded \lsc supersolution of the same problem. For any 
    $x\in\R^N$ and $r>0$, the following result holds:
    \begin{equation}\label{ineq:lcr}
        \max_{B_{r}(x)}(u-v)^{+} \leq \max_{\partial B_{r}(x)}(u-v)^{+}.
    \end{equation}
\end{corollary} 
 
\begin{proof}
    Given $x\in\R^N$ and $r>0$, we apply Proposition~\ref{tcompracao} with $V:=B_r(x)$. It follows
    that $v$ is a supersolution of the Ishii problem associated to $A^{m_V}$, for some $m_V\in N_*$.
    On the other hand, as we already noticed, $u$ is obviously a stratified subsolution of the
    Ishii problem associated to this same compact control set $A^{m_V}$.
    So, applying Theorem~\ref{tcompracaoppp} with $m=m_V$ we deduce that \eqref{ineq:lcr} holds.
\end{proof}

%-----------------------------------------------------------
\subsection{Global comparison results}\label{subsect:gcr}

As explained in \cite{Barles2018}, a (GCR) can be reduced to a (LCR) which
is simpler, by requiring two additional assumptions. The framework in \cite{Barles2018} is far more
general than what we need here and leads to some complexities that are not necessary in the present
situation. For this reason, we provide here a simplified version, well adapted to our needs,
introducing \LOCa and \LOCb below.

Let $\mathcal{C}_{*}$ be a set of subsolutions with certain growing hypotheses, and let $\mathcal{C}^{*}$ be 
a set of supersolutions with certain growing hypotheses. Of course, we have in mind \textit{stratified}
subsolutions here. Here are the localization assumptions we will use:
 
 \

\noindent\LOCa: \textit{Given an \usc subsolution $u \in \mathcal{C}_{*}$ and a \lsc supersolution $v
\in\mathcal{C}^{*}$ there exists a sequence of \usc subsolutions $(u_{\beta})_{\beta}$ such that
for all $x\in\R^N$,}
$$
     \lim\limits_{|x|\to \infty}(u_{\beta}-v )(x)= -\infty\quad\mbox{and}\quad 
     \lim\limits_{\beta\to 0}u_{\beta}(x) = u(x).
$$
 
\noindent\LOCb: \textit{For any $x \in \mathbb{R}^{N}$ and $r>0$, if $u$ is an \usc subsolution, 
    there exists a sequence $(u_{\gamma})_{\gamma>0}$ of
    \usc subsolutions such that for each $\gamma>0$, 
  \[u_{\gamma}(x)-u_{}(x)\geq u_{\gamma}(y)-u_{}(y) + d(\gamma) 
  \quad\mbox{for all } y \in \partial B_r(x)\,,\mbox{ for some } d(\gamma)>0.\]  
  Moreover, $u_{\gamma}(z) \to u_{}(z)$ as $\gamma \to 0$ for any $z \in  B_r(x)$.}
 
\

The growth conditions included in $\mathcal{C}_*$ and $\mathcal{C}^*$ have to be specified in each
case. Notice that in \cite{Barles2018}, the fact that the control spaces are bounded implies
implicitly that the various comparison results are obtained for \textit{bounded} sub and
supersolutions, so that no growth conditions are required. 

But, for instance, in the quadratic example \eqref{ex:quadratic}, the typical growth assumption for
sub and supersolutions is that both should grow strictly less than quadratically. In that case, \LOCa and
\LOCb can be obtained by using 
$$u_\beta(x):=\beta u(x)+(1-\beta)|x|^2\quad\text{and}\quad
 u_{\gamma}(x):=u(x)+\gamma\big(|x-x_0|^2)\,.$$
We refer to Section~\ref{examples} for a detailed and more general example.

For the sake of completeness we provide here the proof that with those assumptions, (LCR) implies
(GCR), this is, subsolutions and supersolutions of the Ishii problem are ordered:

\begin{proposition}\label{prop:lcr.gcr}
    Under   hypotheses \LOCa and \LOCb, a Local Comparison Result implies a Global Comparison
    Result.
 \end{proposition}
 
\begin{proof} 
    Let $u\in\mathcal{C}_*$ and $v\in\mathcal{C}^*$ be respectively an \usc subsolution and a \lsc
    supersolution.  Thanks to \LOCa there exists $x_{\beta}\in\R^N$  such that 
    \begin{equation}\label{prop_4_4_0} 
    M_{\beta}:=\max_{\mathbb{R}^{N}}(u_{\beta}-v)= u_{\beta}(x_{\beta})-v(x_{\beta})\,. 
    \end{equation}
    We assume that $M_\beta>0$ and fix some $r>0$ for which the (LCR) holds around $x_\beta$.
    Notice that by the maximum point property, in particular
    $$M_\beta=u_\beta(x_\beta)-v(x_\beta)\leq \max_{\partial B_r(x_\beta)}(u_\beta-v)\,.$$

    Now thanks to \LOCb, there exists a sequence of subsolutions $(u_{\beta\gamma})_{\gamma}$ that are  
    approximations of $u_{\beta}$ in $B_{r}(x_\beta)$ and using the Local Comparison Result with
    $u_{\beta\gamma}$ and $v$ in $B_r(x_{\beta})$ yields
   \begin{equation}\label{prop_4_4_1}
   \max_{\bar B_r(x_\beta)}(u_{\beta\gamma}-v)^{+} \leq 
    \max_{\partial B_r(x_\beta)}(u_{\beta\gamma}-v)^{+}.
    \end{equation}
    Notice that since $M_{\beta,\gamma}:=\max_{\bar B_r(x_{\beta})}
    (u_{\beta\gamma}-v)\to M_\beta>0$, for $\gamma$ close enough to zero we can leave out
    the positive part of the functions in the previous inequality \eqref{prop_4_4_1}.
    Consequently, thanks to \eqref{prop_4_4_0} and \LOCb,  we have that 
     $$\begin{array}{ll}
         u_{\beta,\gamma}(x_{\beta})-v(x_{\beta}) & \leq \max\limits_{\partial B_r(x_{\beta})}
         (u_{\beta,\gamma}-v)\\
    & \leq  \max\limits_{\partial B_r(x_{\beta})} (u_{\beta}-v) + \max\limits_{\partial
         B_r(x_{\beta})}(u_{\beta,\gamma}-u_{\beta}) \\
    & \leq M_{\beta} + (u_{\beta,\gamma}(x_{\beta})-u_{\beta}(x_{\beta}))- d(\gamma)\,,
    \end{array}$$
    leading to $M_\beta\leq M_\beta-d(\gamma)<M_\beta$ which is a contradiction.
    Therefore, $M_{\beta} \leq 0$ for any $\beta>0$ small enough and taking limits as $\beta\to0$ 
    in \eqref{prop_4_4_0}, we end up with $u(x)\leq v(x)$ for any $x \in \R^{N}$.
 \end{proof}

As a consequence of the previous Proposition~\ref{prop:lcr.gcr}, we obtain the uniqueness of a stratified solution 
of the Ishii problem:
\begin{corollary}\label{cor:comp}
    Assume \HA, \HB, \HC, \LOCa and \LOCb. Then a Global Comparison Result holds between locally
    bounded stratified subsolutions and supersolutions of \eqref{ISHII}. As a consequence, there
    exists a unique stratified solution of \eqref{ISHII}.
\end{corollary}

%------------------------------------------
\subsection{Extremal Ishii Solution}

We prove here that $U_A^-$ and $U_A^+$ are respectively the minimal and maximal solutions in the
sense of Ishii. Of course, we recall that $U_A^-$ is the unique \textit{stratified} solution, the
only one enjoying the complementary inequality $H_T\leq0$ among all Ishii subsolutions.

\begin{proposition}\label{U-comparasion}
    Under hypotheses \HA, \HB, \HC, \LOCa, \LOCb,\\[2mm] 
    $(i)$ $U_A^{-}$ is the minimal
    locally bounded \usc supersolution (and solution) of the Ishii problem \eqref{ISHII}; \\[2mm]
    $(ii)$ $U_A^{+}$ is the maximal
    locally bounded \lsc subsolution (and solution) of the Ishii problem \eqref{ISHII};\\[2mm]
    $(iii)$ $U_A^-$ is the unique stratified solution of \eqref{ISHII}.
\end{proposition}

\begin{proof}
    Let $v$ be any supersolution  of \eqref{ISHII}. Since $U_A^-$ is a stratified subsolution, we can
    apply the global comparison result to $U_A^-$ and $v$, then we have that $U_A^-\leq v$ in $\R^N$. 
    As $U_A^-$ is itself of course a supersolution, it
    is clearly the minimal one. Of course, since any viscosity solution is a supersolution by
    definition, then for any solution $u$, the argument proves that $U_A^-\leq u$ so that $U_A^-$ is
    also the minimal viscosity solution.

    Concerning the maximal subsolution, the simplest way to obtain the result it is to recall that
    if $u$ is any viscosity subsolution in the sense of Ishii (we do not require here that it is a
    stratified subsolution), it is also a viscosity subsolution of the Ishii problem associated to
    the compact control spaces $A^m$ for any $m\in\N_*$.

    Using the fact that $U_{A^m}^+$ is the maximal subsolution for the bounded control case (see
    \cite[Corollary 4.4]{manoel}), we deduce that $u\leq U_{A^m}^+$ for any $m\in\N_*$. Then, passing to the limit
    as $m\to+\infty$ yields the result : $u\leq U_A^+$.

    The proof of $(iii)$ just follows from Corollary~\ref{cor:comp}. Since $U_A^-$ is a stratified
    solution of \eqref{ISHII}, it is the unique one.
\end{proof}

%%%%%%%%%%%%%%%%%%%%%%%%%%%%%%%%%%%%%%%%%%%%%%%%%%%%%%%%%%%%%%%%%%%%%%%%%%%
\section{A continuous family of solutions}\label{sect:U.eta}

%%%%%%%%%%%%%%%%%%%%%%%%%%%%%%%%%%%%%%%%%%%%%%%%%%%%%%%%%%%%%%%%%%%%%%%%%%

In this section we build a whole family of value functions, $U^{\eta}$, with $\eta>0$. Those turn
out to be locally Lipschitz Ishii solutions which yield a continuous path between $U_A^{-}$ and 
$U_A^{+}$. 

To build them, let us start with defining $\eta$-trajectories. For $\eta>0$ we set 
\begin{equation}\label{-+controleta1} 
    \begin{aligned}
        \tau^{\eta}_{A}(x):=\Big\{  (X(\cdot),\alpha(\cdot)) \in & \tau_A(x)~|~ \mbox{for a.e.}
         \ t\in \mathcal{E}_{\mathcal{H}},\\ & b_2(X(t),\alpha(t))\cdot e_N \geq -\eta,
         ~ b_1(X(t),\alpha(t))\cdot e_N \leq \eta \Big\}.
    \end{aligned}
\end{equation}
All $\eta-$trajectories are of course not necessarily regular, however they are close to regular for
$\eta>0$ small. Notice that for all $0<\eta\leq\eta'$,
$$\tau_A^{\mathrm{reg}}(x)\subset \tau^{\eta}_{A}(x)\subset\tau_A^{\eta'}(x)\subset\tau_A(s)\,.$$ 
Now, let us define the associated $\eta$-value function as usual:
\begin{equation}\label{ueta}
    U^{\eta}_{A}(x) = \inf_{\tau^{\eta}_{A}(x)}\left( \int_{0}^{\infty}
    l( X_x(t),\alpha(t)) e^{-\lambda t} dt \right).
\end{equation}
Of course, for $0<\eta\leq \eta'$ as above,
\begin{equation}\label{valuesfunctions}
U_{A}^{-} \leq U_{A}^{\eta'}\leq U_{A}^{\eta}  \leq U_{A}^{+}.
\end{equation} 

Having a look at Proposition~\ref{continuouss}, it is easy to see that the same arguments work for any
$\eta$-value function. So, we claim that $U^\eta_A$ is locally bounded and locally Lipschitz continuous.
We omit the details here since we actually prove a stronger result below, see
Lemma~\ref{lem:uni.lip}: the $(U^\eta_A)_\eta$ are locally uniformly Lipschitz and bounded.

%------------------------------------------------------------
\subsection{$\eta$-value functions are Ishii solutions}

As is standard, let us begin with the Dynamic Programming Principle:
\begin{theorem}[Dynamic Programming Principle]\label{dppeta} 
    Let $\eta>0$. Then for any $T>0$, 
    \begin{equation}\label{prdineta0}
        U_{A}^{\eta}(x)  = \inf_{\tau^{\eta}_{A}(x)}\left( \int_{0}^{T} 
        l( X(t),\alpha(t)) e^{-\lambda t} dt + U_{A}^{\eta} \left(X\left(T \right) 
        \right)e^{-\lambda T}\right). 
    \end{equation}
\end{theorem}

\begin{proof} 
    The proof is analogous to Theorem \ref{dppt} in section \ref{sect:VF}.
\end{proof}

We also define a $\eta$-tangential Hamiltonian that will play the role of $H_T$ or
$H_T^\mathrm{reg}$ for $U^\eta_A$. Let 
$$
    A_{0}^{\eta}(x)=\{\alpha=\left(\alpha_1,\alpha_2, \mu \right) \in  A_0 
    ~|~b_1(x,\alpha_1,\alpha_2)\cdot e_N\le \eta,  ~b_2(x,\alpha_1,\alpha_2)\cdot e_N\ge -\eta\}\,,
$$ 
we define the $\eta$-tangential Hamiltonian as
\begin{equation}\label{eta_tanH}
    H^{\eta}_{T}(x,\phi(x), D_{\mathcal{H}}\phi)=
    \sup\limits_{\alpha\in A^{\eta}_0(x)}   \big\{
    \lambda \phi \left(  x \right)  - b_{\mathcal{H}}
    \left(x, \alpha\right)\cdot  \left( D_{\mathcal{H}}\phi \left( x\right),0\right)-
    l_{\mathcal{H}}( x,\alpha)\big\}
\end{equation}
and notice that Proposition~\ref{reducao} obviously applies to $H_T^\eta$ as well.

Let us prove below that the $\eta$-value functions are Ishii solutions. 
\begin{theorem}\label{existenciaeta} 
    Assume that \HA, \HB, \HC hold and let $\eta>0$. Then the value function $U_A^{\eta}$ is a
    viscosity solution of the Ishii problem, \eqref{ISHII}. Moreover, $U^\eta_A$ satisfies
    $H^{\eta}_{T}(x,u,Du)\leq0$ on $\H$.
\end{theorem}

\begin{proof}
    Since the Dynamic Programming Principle \ref{dppeta} holds, the proof that $U_A^\eta$ is at the
    same time a sub and and supersolutions of \eqref{ISHII} is analogous to the standard case, see
    Theorem~\ref{existencia}.

    In order to prove the $H^\eta_T$ property we need some adaptations which follow the lines of
    Theorem~\ref{existencia}-$(ii)$. Here we face several cases as follows:
    \begin{enumerate}
      \item $b_1(x,\alpha_1)\cdot e_N < \eta \quad \mbox{and} \quad b_2(x,\alpha_2)\cdot e_N > -\eta.$
      \item $b_1(x,\alpha_1)\cdot e_N =\eta\quad\mbox{and}\quad b_2(x,\alpha_2)\cdot e_N = -\eta.$
      \item	$ b_1(x,\alpha_1)\cdot e_N = \eta\quad \mbox{and}\quad b_2(x,\alpha_2)\cdot e_N > -\eta.$ 
      \item	$ b_1(x,\alpha_1)\cdot e_N < \eta\quad \mbox{and}\quad b_2(x,\alpha_2)\cdot e_N =-\eta.$
    \end{enumerate}
    We do not detail more the computations since they are obvious adaptations. Here also we use the
    controllability for dealing with the limit cases $b_i\cdot e_N=\eta$ or $-\eta$. The reader can
    also check similar proofs in \cite[Theorem 2.5]{manoel} for the case $\eta=0$.
\end{proof}

%------------------------------------------------------------
\subsection{Asymptotics as $\eta\to+\infty$ of $U_A^\eta$}

In this section we prove the asymptotic results for the Hamiltonian $H_T^\eta$ and the $\eta$ value function $U_A^\eta$. 

Let us begin with the results related to the Hamiltonian $H_T^\eta$.

\begin{lemma}\label{conv.ht}
    Assume \HA, \HB and \HC. Then the following limits are monotone and locally uniform in
    $\R^N\times\R\times\R^N$: 
    $$\lim_{\eta\to+\infty}H_T^\eta=H_T\,,\qquad\lim_{\eta\to0^+}H_T^\eta=H_T^\mathrm{reg}\,.$$
\end{lemma}
\begin{proof}
    It is clear from the definition of $A_0^\eta$ that this sequence of sets is monotone with respect
    to $\eta$ and that moreover $A_0^\eta\to A_0$ as $\eta\to+\infty$ while $A_0^\eta\to
    A_0^\mathrm{reg}$ as $\eta\to0$.
    On the other hand, from Proposition~\ref{reducao}, we know that on any fixed compact
    $K\subset\R^N\times\R\times\R^N$, the sup defining $H_T^\eta$ is attained on a compact control
    set $A_K\subset A_0$ which can be chosen uniformly with respect to $\eta$. 

    \smallskip

    \indent Let us consider the case $\eta\to+\infty$. We already know that obviously
    $H_T^\eta\leq H_T$ since the control sets involved are ordered the same way. Now let us fix a
    compact set $K\subset\R^N\times\R\times\R^N$. There exists $\tilde A_K\subset A_0$ such that for any 
    $(x,u,p)\in K$, we have an optimal control $a\in \tilde A_K$. In other words,
    $$H_T(x,u,p)=\lambda u-b_\H(x,a)\cdot p-l_\H(x,a)\,.$$
    But since $\tilde A_K$ is compact, it follows that there exists
    $\eta=\eta(K)$ big enough such that $\tilde A_K\subset A_0^{\eta}$. In other words,
    for any $(x,u,p)\in K$, $H_T(x,u,p)\leq H_T^\eta(x,u,p)$ since this last Hamiltonian is taken as
    the supremum over $A_0^\eta$. This implies that in fact $H_T=H_T^\eta$ on $K$ and of course the
    same property holds for $\eta'>\eta$. The asymptotic result follows.

    \smallskip

    \indent We turn now to the case $\eta\to0$. Here the inequality $H_T^\eta\geq
    H_T^\mathrm{reg}$ is the obvious one. To get the reverse, we use a similar approach as above,
    but with some adaptations.

    Since $H_T^\eta$ is non-increasing with respect to $\eta$ and bounded from below by
    $H_T^\mathrm{reg}$, the following limit is well-defined for any $(x,u,p)$
    $$h(x,u,p):=\lim_{\eta\to0}H_T^\eta(x,u,p)\,.$$ 
    Now consider $K$ compact as above and $(x,u,p)\in K$. Given a decreasing sequence $\eta_n\to0$, 
    there exists a sequence of associated optimal controls $(a_n)_n$ for $H_T^{\eta_n}(x,u,p)$.
    As we already did above, there exists a compact set $\tilde A_K\subset A_0^{\eta_n}$ uniformly with respect
    to $n$, such that $a_n\in \tilde A_K$ for any $n\in\N$.

    So, we can extract a subsequence still denoted by $(a_n)_n$ converging to some control $a_*\in
    \tilde A_K$. Since for any $n$, the optimal control $a_n\in A_0^{\eta_n}$, we deduce that $a_*\in A_0^\mathrm{reg}$ and
    passing to the limit in the Hamiltonian we get
    $$\lim_{n\to+\infty} H^{\eta_n}(x,u,p) = \lambda u-b_\H(x,a_*)\cdot p-l_\H(x,a_*)\leq
    H_T^\mathrm{reg}(x,u,p)\,.$$
    This implies in the end that $\lim_{\eta\to0} H_T^\eta = H_T^\mathrm{reg}$, and the limit is uniform on $K$
    because the controls involved remain in a fixed compact set, and $(b,l)$ are continuous
    with respect to $a$.
\end{proof}

\medskip

\begin{remark}
    Actually we have proved a stronger property for the case $\eta\to+\infty$: for any compact
    $K\subset\R^N\times\R\times\R^N$ there exists $\eta_0(K)>0$ such that
    $$\text{for any }\eta\geq\eta_0\,,\quad H_T^\eta=H_T\text{ on }K\,.$$
\end{remark}

\medskip

\begin{lemma}\label{lem:uni.lip}
    Assume that \HA, \HB and \HC hold.
    The sequence $(U_A^\eta)_{\eta>0}$ is non-increasing with respect to $\eta$, locally uniformly
    bounded and locally uniformly Lipschitz in $\R^N$ with respect to $\eta$. As a consequence, the
    following limits
    \begin{equation}\label{def:u}
        w^-(x):=\lim_{\eta\to+\infty}U_A^\eta(x)\quad\text{ and }\quad
        w^+(x):=\lim_{\eta\to0}U_A^\eta(x)
    \end{equation}
    are classical (Ishii) viscosity solutions of \eqref{ISHII}.
\end{lemma}
\begin{proof}
    Let us first notice that the monotonicity property is obvious since the sequence of controlled
    trajectories $\tau^\eta_A$ is non-decreasing as $\eta$ increases. As we already noticed,
    for any $\eta>0$, $U_A^-\leq U_A^\eta\leq U_A^+$, so the sequence is locally
    uniformly bounded in $\R^N$.

    Taking a look at the proof of Proposition~\ref{continuouss}, it is clear that the same
    computations are valid for $U_A^\eta$, and the local Lipschitz constant $C(V)$ does not dependent on
    $\eta$ since it only depends on the local bound of $U_A^\eta$, and the parameter
    $\delta>0$ in \HB. So, the result holds.

    Thanks to the former properties, the two functions $w^-$ and $w^+$ are well-defined and
    continuous functions in $\R^N$. Moreover, by the standard stability properties of viscosity
    solutions, it is clear that both are viscosity solutions in the sense of Ishii of \eqref{ISHII}.
\end{proof}

Now we can prove easily one of the convergence result :

\begin{proposition}
    Assume \HA, \HB, \HC, \LOCa, \LOCb. Then the following limit holds
    $$ \lim_{\eta\to+\infty}U_A^\eta=U_A^-\quad\text{locally uniformly in }\R^N\,.$$
\end{proposition}
\begin{proof}
    By Lemma~\ref{lem:uni.lip}, we already know that $w^-$ is a viscosity solution in the sense of
    Ishii of \eqref{ISHII}. Moreover, $H_T^\eta$ converges locally uniformly to $H_T$, we also
    deduce that $w^-$ satisfies the subsolution inequality $H_T(x,w^-,Dw^-)\leq0$ in the viscosity
    sense on $\H$. In other words, $w^-$ is a \textit{stratified} solution of \eqref{ISHII}. But by
    uniqueness of such solution (see Proposition~\ref{U-comparasion}-$(iii)$), it follows that
    $w^-\equiv U_A^-$ and the result is proved.  
\end{proof}

\begin{remark}
    Another proof can be obtained by performing a detailed analysis of the trajectories involved in
    the definition of the values functions, proving that $U_A^\eta\to U_A^-$ directly without using
    the comparison argument. However, this analysis requires some extra assumptions on the elements
    $(b,l)$ in order to estimate the spreading of optimal trajectories. This is the approach that
    we use below when $\eta\to0$. Here, on the contrary, using the comparison argument for stratified
    solutions yields an easy proof. Notice however that we need hypotheses \LOCa and \LOCb to
    use such arguments.
\end{remark}

%------------------------------------------------------------
\subsection{Asymptotics as $\eta\to0$ of $U_A^\eta$}

Unfortunately, getting a similar result when $\eta\to0$ is not so easy since it requires to check
that the convergence of $\eta$-trajectories yield regular trajectories as $\eta\to0$. What seems
like an obvious result actually requires a fine proof, adapting  the strategy used in \cite[Lemma
5.3]{manoel2014}. We prove this result for compact control sets $A^m$. To simplify the notation of the controls 
in the lemma below, we do not use the notation $a^m$, and substitute it by $a$,  although the controls belong to $A^m$.

\begin{lemma}\label{lem:eta.reg}
    Assume \HA,  \HB, \HC. Let $\eta_n\to0$ and $(X^n,\alpha^n)\in\tau_{A^m}^{\eta_n}$ be a sequence of $\eta$-trajectories
    defined on $[0,T]$. We assume that $(X^n)$ converges uniformly to some admissible trajectory $X$
    on $[0,T]$. Then the trajectory $X$ is regular. In other words, there exists a control, $a$, such
    that $(X,a)\in \tau_{A^m}^\mathrm{reg}$.
\end{lemma}

\newcommand{\dist}{\mathop{\rm dist}}

\begin{proof}
    Since $X$ is an admissible trajectory, there exists a control $\alpha$ such that
    $(X,\alpha)\in\tau_{A^m}$. Let $z \in \mathcal{H}$, and define
    \[
    \begin{array}{l}
	K(z):=\lbrace b_{\mathcal{H}}(z,\alpha'):  \alpha' \in A_0^{\mathrm{reg},m}(z)\rbrace.
    \smallskip\\
	E_{sing}^{\Upsilon}:=\lbrace s \in [0,T]: X(s)\in \mathcal{H} \mbox{ and } \dist(b_{\mathcal{H}}(X(s),	
	\alpha(s)),K(X(s)))\geq \Upsilon  \rbrace.
    \smallskip\\
        E_{sing}:=\lbrace s \in [0,T]: X(s)\in \mathcal{H} \mbox{ and } 
        \dist(b_{\mathcal{H}}(X(s),\alpha(s)),K(X(s)))>0  \rbrace,
    \end{array}
    \]
    where $\dist(\cdot)$ is the euclidian distance in $\R^N$. 
    Our aim is to prove that $|E_{sing}|=0$, which will prove that the trajectory is regular. To
    this end, since $E_{sing}= \cup_{j \in \mathbb{N}} E_{sing}^{1/j}\,,$ it is enough to prove that
    for any $\Upsilon>0$, $|E_{sing}^{\Upsilon}|=0$.

    Assuming that there exists $\Upsilon >0$ such that $|E_{sing}^{\Upsilon}|>0$, we 
    prove that we can find a control $\tilde{\alpha}(\cdot)\in
    L^{\infty}(E_{sing}^\Upsilon;A^m)$ satisfying the following conditions: for any $s$ in
    $E_{sing}^\Upsilon$, $\tilde{a}(s) \in A_0^{\mathrm{reg},m}(X(s))$ and
    \begin{equation}\label{appppammçç}
         b_{\mathcal{H}}( X^{n}(s),\alpha^{n}(s))=
         b_{\mathcal{H}}( X(s),\tilde{a}(s)) + \varsigma^n(s)
     \end{equation} 
     where $\varsigma^n(\cdot)$ is measurable and goes uniformly to $0$ when $n$ goes to infinity.
     By doing so, we prove that we can redefine the control $\alpha$ so that
     $|E_{sing}^\Upsilon|=0$. Indeed, by passing to the limit as $n\to+\infty$ we see that the
     dynamics associated to $\alpha$ and $\tilde a$ are the same on $E_{sing}^\Upsilon$, so that the
     trajectory $X$ can be also constructed by using the control $\tilde a$, which is regular. 

     Writing $\alpha^n=(\alpha^n_1,\alpha^n_2,\mu^n)$, notice that in what follows, we will not
     modify the parameter $\mu^n$ , only the controls $\alpha_i^n$, and that of course, no change is
     needed if $X_n\notin\H$.

     \smallskip

     \noindent\textbf{Step 1:} The first step consists in finding for each $n$ a \textit{regular} control
     $a^{n}(\cdot)$ for the trajectory $X^{n}(\cdot)$ satisfying on $\Delta_n:=(X_n)^{-1}(\H)\cap
     E_{sing}^\Upsilon$
     $$b_{\mathcal{H}}(X^{n}(s),\alpha^{n}(s)) = b_{\mathcal{H}}(X^{n}(s),a^{n}(s))+o_n(1)\,.$$
     Denoting by $\alpha^n=(\alpha_1^n,\alpha_2^n,\mu^n)$ we set
    \begin{equation}\label{maxmim}
	    \gamma^{n}_1(s):= -\max\lbrace 0, b_1( X^{n}(s),\alpha_1^{n}(s)) \cdot e_N \rbrace, 
	    \qquad  
	    \gamma^{n}_2(s):= -\min\lbrace 0, b_2( X^{n}(s),\alpha_2^{n}(s))\cdot e_N \rbrace
    \end{equation} 
    and we claim that
     \begin{equation}\label{3329}
         \mu^{n} \gamma^{n}_1(s) + (1-\mu^{n} )\gamma^{n}_2(s)=0,
         \mbox{ for a.e.} ~~s\in \Delta_n.
     \end{equation}
    Indeed, if $\alpha^{n}(s)\in A_0^{\mathrm{reg},m}(X^{n}(s))$, then
    $\gamma_1^{n}(s)=0=\gamma_2^{n}(s)$ and the claim is obvious. 
    So, assume that $\alpha^{n}(s) \notin A_0^{\mathrm{reg},m}(X^{n}(s))$. By Theorem \ref{dynamic},
    since $X^{n}(s)\in\mathcal{H}$, we have that $b_{\mathcal{H}}(X^{n}_x(s),\alpha^{n}(s))\cdot
    e_N=0$ almost everywhere on $\Delta_n$ and we face three cases which lead to \eqref{3329}:
    \begin{enumerate}
        \item $\gamma_i^{n}(s):=-b_i( X^{n}_x(s),\alpha_i^{n}(s)) \cdot e_N$~~ if  $\mu^{n}(s) \in (0,1)$;
        \item $\gamma_1^{n}(s)=0$~~ if $\mu^{{n}}(s)=1$;
        \item $\gamma_2^{n}(s)=0$~~ if $\mu^{{n}}(s)=0.$
    \end{enumerate} 
    Notice that since $b_i$ and $X^{n}$ are continuous and $\alpha_i^{n}(\cdot)$ is measurable, then
    $\gamma^{n}_i$ is measurable. Moreover, since $|\gamma^{n}_i(s)|\leq \eta_n$ almost everywhere
    on $\Delta_n$, it follows that
    \begin{equation}\label{gammaitocero}
        \gamma_i^n\to 0,\quad \mbox{ uniformly as }~ n\to\infty.
    \end{equation}
    Furthermore, by construction we have approximating regular dynamics
    \begin{equation}\label{0000}
        (b_1( X^{n}(s),\alpha_1^{n}(s))+ \gamma^{n}_1(s)e_N )\cdot e_N \leq 0, 
        \ \ (b_2( X^{n}(s),\alpha_2^{n}(s))+\gamma^{n}_2(s)e_N) \cdot e_N \geq 0.
    \end{equation}
    Now, we use $\gamma_i^{n}$ to build the regular control $a^{n}(\cdot)$. To do so, we define 
    $$\beta^{n}(s):=\min\left\{ \frac{\delta - 2|\gamma^{n}_1(s)| }{\delta}, 
    \frac{\delta - 2|\gamma^{n}_2(s)| }{\delta} \right\}.$$ 
    Thanks to \eqref{gammaitocero}, we have that $\beta^{n}(s) \to 1$ when $n \to \infty$. 

    \smallskip

    \noindent\textbf{Step 2:} We claim that there exists $\tilde\alpha_i^{n}(s) \in A^m_i$ 
    such that
    \begin{equation}\label{conbbvvvv}
        b_i( X^{n}(s),\alpha_i^{n}(s))= b_i( X^{n}(s),\tilde\alpha_i^{n}(s)) + p_i^{n}(s),
    \end{equation}
    where $p_i^{n}(s):=\left(1-\beta^n(s)\right)
    \left(b_i( X^{n}(s),\alpha_i^{n}(s))+\gamma^{n}_i(s) \right) - \gamma^{n}_i(s).$ 
    Note that $p^{n}_i$ is measurable and goes uniformly to $0$ when $n$ goes to infinity.
    Indeed, if $\beta^{n}(s)=1$, then $|\gamma^{n}_1(s)|=0=|\gamma^{n}_2(s)|$. In this case we can take 
    $$\tilde\alpha_i^{n}(s):=\alpha_i^{n}(s)~~ \mbox{and}~~
     p_i^{n}(s)=(1-1)(b_i( X^{n}_x(s),\alpha_i^{n}(s))+ 0 ) - 0=0.$$ 

    Otherwise, assume that $\beta^{n}(s)\neq 1$. 
    Since $\kappa_n(s):=\beta^{n}(s)\frac{\gamma^{n}_i(s) e_N}{1-\beta^{n}(s)}\in\{z\in\R^N:
    |z|\leq\frac{\delta}{2}\}$, for $n$ big enough,  it follows that
    $$ \beta^{n}(s) \left( b_i( X^{n}(s),\alpha_i^{n}(s))+ \gamma^{n}_i(s)e_N \right)=
    \beta^{n}(s)b_i( X^{n}(s),\alpha_i^{n}(s))+ (1-\beta^{n}(s))\kappa_n(s)$$
    belongs to $\mathcal{B}(X^n(s))$. This is due to \HB and the convexity property of the images
    $\mathcal{B}(X(s))$.
    In other words, there exists $\tilde\alpha_i^{n}(s)\in A_i^m$ such that
    $$
    \beta^{n}(s) \left( b_i( X^{n}(s),\alpha_i^{n}(s))+ \gamma^{n}_i(s)e_N \right)=
     b_i( X^{n}(s),\tilde\alpha_i^{n}(s)). $$
    
    Furthermore, $\tilde\alpha^n$ is a regular control.
    Indeed, thanks to \eqref{0000}
    $$b_1( X^{n}(s),\tilde\alpha_1^{n}(s))\cdot e_N= 
    \beta^{n}(s)(b_1( X^{n}(s),\alpha_1^{n}(s))+\gamma^{n}_1(s) e_N)\cdot e_N \leq 0 $$
    $$b_2( X^{n}(s),\tilde\alpha_2^{n}(s))\cdot e_N= \beta^{n}(s)(b_2( X^{n}(s),\alpha_2^{n}(s))+ 
    \gamma^{n}_2(s) e_N)\cdot e_N \geq 0 $$
    and thanks to Theorem \ref{dynamic} and \eqref{0000} 
    \[
    \begin{array}{l}
	\mu^{n}(s)b_1( X^{n}(s),\tilde\alpha_1^{n}(s))\cdot e_N + (1-\mu^{n}(s))b_2( X^{n}(s),
	\tilde\alpha_2^{n}(s))\cdot e_N  
    \smallskip\\
	=\mu^{n}(s)\beta^{n}(s)b_1( X^{n}(s),\alpha_1^{n}(s))\cdot e_N  +(1-\mu^{n}(s)) \beta^{n}(s)
	(b_2(X^{n}(s),\alpha_2^{n}(s))
    \smallskip\\
	~~+ \beta^{n}(s) \mu^{n} (s)\gamma^{n}_1(s) +  \beta^{n}(s)(1-\mu^{n}(s) )\gamma^{n}_2(s)
	\smallskip\\
        =0 \text{  a.e on $\Delta_n$.}
    \end{array}
    \]

     Notice that we constructed
     $(\tilde\alpha_1^{n}(\cdot),\tilde\alpha_2^{n}(\cdot))$ pointwise for $s\in
     \Delta_n$, so that they may not necessarily be measurable. However, thanks to a
     measurable selection argument (Filippov's Lemma) we can find measurable controls
    $a^{n}_i(\cdot)$ such that 
    $$b_i( X^{n}(s),\tilde\alpha_i^{n}(s))-p_i^{n}(s)=b_i( X^{n}(s),a^{n}_i(s)). $$ 
    Thus, $a^{n}(\cdot):=(a^{n}_1(\cdot), a^{n}_2(\cdot), \mu^n(\cdot)) \in A_0^{\mathrm{reg},m}(X^{n}(s))$ 
    is measurable. Consequently,
    $$ b_{\mathcal{H}}( X^{n}(s),\tilde\alpha^n(s)) - p^{n}(s)= b_{\mathcal{H}}( X^{n}(s),a^{n}(s)) $$
    where $p^{n}(s)=\mu^{n}(s)p_1^{n}(s) + (1-\mu^{n}(s))p_2^{n}(s)$. 
    Note that $p^{n}$ is measurable and goes uniformly to $0$ when $n$ goes to infinity 
    since $p_i^n$ goes to zero as $n$ goes to infinity.

    \smallskip

     \noindent\textbf{Step 3:} Since $(a^{n}_1(s), a^{n}_2(s), \mu^{n}(s)) \in
     A_0^{\mathrm{reg},m}(X^{n}(s))$, thanks to \cite[Lemma 5.3]{manoel2014} there exists a measurable
     control $\tilde{a} \in A_0^{\mathrm{reg},m}(X(s))$, i.e. a control associated to the limit
     trajectory, satisfying
     \begin{equation}\label{apppo}
        b_{\mathcal{H}}( X^{n}(s),a^{n}(s))= b_{\mathcal{H}}( X(s),\tilde{a}(s))+ \sigma^{n}(s).
    \end{equation}
    with $\sigma^n(\cdot)$ measurable that goes uniformly to $0$ when $n$ goes to infinity. 
    Therefore,
    $$ b_{\mathcal{H}}( X^{n}(s),\tilde\alpha^{n}(s))= b_{\mathcal{H}}( X(s),\tilde{a}(s))+
    \sigma^{n}(s)+p^{n}(s)\,,$$ 
    leading to \eqref{appppammçç} with $\varsigma^n(s):=\sigma^{n}(s)+p^{n}(s)$.
    As we said, this modification allows to associate the trajectory $X$ with a regular control on
    $E_{sing}^\Upsilon$. By doing the same procedure on each $E_{sing}^{1/j}$, we end up with the
    fact that for some control $a$ such that $(X,a)\in\tau_{A^m}$, $|E_{sing}|=0$, the trajectory $(X,a)$ is regular.
\end{proof}

Next we need some assumptions on the cost and dynamics so that we can control the asymptotic
behaviour of trajectories.

\begin{lemma}\label{TTTTT} 
Assume \HA,  \HB, \HC and let $m\in\N_*$. Assume there exist $\lambda, C$ and $\overline{C}$ such that
\begin{enumerate}
	 \item[\HD] $|b_i(x,\alpha_i^m)|\le \lambda(1+|\alpha_i^m|+|x|)$ and 
         $l_i(x,\alpha_i^m)\leq C|x|^{1-\epsilon} + {\overline C}|\alpha_i^m|$,  
         $\forall x\in\R^N$, $\alpha_i^m\in A_i^m$, $i=1,2$. 
\end{enumerate}
    Let $(X,\alpha) \in \tau_{A^m}(x)$. Then,
\[
	\lim\limits_{T \to \infty}U_{A^m}^{\eta} (X(T))e^{-\lambda T}
	=\lim\limits_{T \to \infty}U_{A^m}^{-}(X(T))e^{-\lambda T}=\lim\limits_{T \to \infty}U_{A^m}^{+}(X(T))e^{-\lambda T}=0.
\]
\end{lemma}
\begin{proof}
    Denote $U_{A^m}^{+}$, $U_{A^m}^{-}$ and $U_{A^m}^{\eta}$ by $U$ for simplicity since the argument is
    the same for all of them. Actually, it is enough prove it for $U_{A^m}^+$ only. 

    Thanks to hypotheses \HB and  \HC, there exists a regular control
    $\alpha^*=(\alpha_1^*,\alpha_2^*,\mu^*)\in A^m$ such that 
    $b_i(X(T),\alpha_i^{*})=0$. The associated trajectory is just $X(t)=X(T)$ for all $t\geq T$.
    By definition of $U$ as the infimum and hypothesis \HD we get
    $$
        U(X(T))e^{-\lambda T} \leq \frac{l(X(T),\alpha_{*})}{\lambda}e^{-\lambda T} 
        \leq  \frac{C|X(T)|^{1-\epsilon} + \overline{C}|\alpha_{*}|}{\lambda}e^{-\lambda T}.
    $$
    Thanks to \HD, it follows from Gronwall's lemma on $X(\cdot)$ that $|X(T)|\leq |x+c| e^{\lambda T}$ so that 
    $$U(X(T))e^{-\lambda T}\leq \frac{1}{\lambda}\Big(C|x+c|^{1-\epsilon}e^{\lambda(1-\epsilon)
    T}+\bar C|\alpha_*|\Big)e^{-\lambda T}\to0\quad\text{as }T\to+\infty\,,$$
    which proves the result.
\end{proof}

The following theorem states the convergence of $U_{A^{m}}^{\eta}$ to $U_{A^{m}}^{+}$ as $\eta$ goes
to zero on compact control sets. 
\begin{theorem}\label{cooo} 
	Under the hypotheses  \HA, \HB, \HC and \HD, it follows that
	$$\lim_{\eta\rightarrow 0}U_{A^m}^{\eta} (x)=U_{A^m}^{+}(x).$$
\end{theorem}
\begin{proof}
    Let us consider a strictly positive decreasing sequence $\lbrace \eta_n \rbrace_{n\ge 0}$
    such that $\lim\limits_{n\to \infty} \eta_n=0$. We assume that for some $x\in\R^N$,
    \begin{equation}\label{UetamenorUmais}
        U_{A^m}^{\eta_n}(x)< U_{A^m}^{+}(x),
    \end{equation}
    otherwise, there exists $N\in \N$ such that $U_{A^m}^{\eta_n}(x)= U_{A^m}^{+}(x)$ for all $n\ge
    N$ and the result holds.

    \smallskip

    \noindent\textbf{Step 1: }
    Given $T>0$, thanks to \eqref{UetamenorUmais} and   Theorem \ref{dppeta}, there exists
    $(X_{x}^{\eta_n},\alpha^{\eta_n}) \in \tau^{\eta_n}_{A^m}(x)$ satisfying 
    \begin{equation}\label{nnnnetatt}
        \int_{0}^{T} l( X^{\eta_n}_x(t),\alpha^{\eta_n}(t)) e^{-\lambda t} dt + 
        U_{A^m}^{\eta_n}( X^{\eta_n}_x(T))e^{-\lambda T}\leq U_{A^m}^{\eta_n}(x) +
        \frac{U_{A^m}^{+}(x)-  U_{A^m}^{\eta_n}(x)}{2} < U_{A^m}^{+}(x).
    \end{equation}
    Let us define 
    $$Y^{\eta_n}(s):=\int_{0}^{s}l(X_x^{\eta_n}(t),\alpha^{\eta_n}(t))dt.$$  
    Since $b_i(X_x^{\eta_n}, \alpha^{\eta_n})$ and $l_i(X_x^{\eta_n}, \alpha^{\eta_n})$ are bounded
    in $[0, T]$, it follows that the curves $(X_x^{\eta_n},Y^{\eta_n})(\cdot)$ are equicontinuous
    and uniformly bounded in $[0, T]$. Therefore, by Ascoli Arzela's Theorem, we can extract a
    subsequence $(X_x^{\eta_n}, Y^{\eta_n})(\cdot)$ which converges uniformly to
    $Z^{T}:=(X_x^{T},Y^{T})$ in $[0, T]$. Proceeding as in \cite[Lemma 5.3]{manoel2014}, we find a
    measurable control $\alpha^{T}(\cdot)$ such that
    \begin{equation}
         (\dot{X}_x^{T}(t), \dot{Y}^{T}(t))= 
         (b(X_x^{T}(t),\alpha^{T}(t)),l(X_x^{T}(t),\alpha^{T}(t)))\quad\forall t\in[0,T].
    \end{equation}
    Thus,  $\left((b(X^{\eta_n}(t), \alpha^{\eta_n}(t)),l
    (X^{\eta_n}(t),\alpha^{\eta_n}(t))\right))$ converges weakly$^*$ to
    $\left(b(X^{T}(t),\alpha^{T}(t)),l(X^{T}(t),\alpha^{T}(t))\right)$ in
    $L^{\infty}([0,T];\mathbb{R}^{N+1})$. Moreover, by Lemma~\ref{lem:eta.reg}, we may assume that
    $\alpha^T$ is regular.

    \smallskip

    \noindent\textbf{Step 2:} Since $l (X_x^{\eta_n}(t),\alpha^{\eta_n}(t))$ converges weakly$^*$ to
    $l(X_x^{T}(t),\alpha^{T}(t))$ in $L^{\infty}([0,T];\mathbb{R})$ and the exponential  is in
    $L^{1}([0,T];\mathbb{R})$, we obtain that
    $$\lim_{n\rightarrow \infty}  \int_{0}^{T} l( X^{\eta_n}_x(t),\alpha^{\eta_n}(t))
    e^{-\lambda t} dt =  \int_{0}^{T} l( X^{T}_x(t),\alpha^{T}(t)) e^{-\lambda t} dt.$$ 
    Thanks to hypotheses \HB there exists a bounded control $\alpha_i^{*}$ such that
    $b_i(X_x^{T}(T),\alpha_i^{*})=0 $, which allows to consider the trajectory
    $X_x^{T}(t):=X_x^{T}(T)$ associated to $\alpha^{T}_i(t):=\alpha_i^{*}$ for any $t\geq T$. 
    Thus, $( X^{T}_x(\cdot),\alpha^{T}(\cdot))
    \in \tau_{A^m}^{\mathrm{reg}}(x)$ and thanks to the Dynamic Programming Principle,
    \eqref{UetamenorUmais} and \eqref{nnnnetatt},  
    \begin{equation}\label{inUmais11}
        \begin{aligned}
            U_{A^{m}}^{+}(x) -  U_{A^{m}}^{+}( X^{T}_x(T))e^{-\lambda T}  
            &\leq \int_{0}^{T} l( X^{T}_x(t),\alpha^{T}(t)) e^{-\lambda t} dt \\
            &= \lim\limits_{n\rightarrow \infty}\int_{0}^{T} 
            l( X^{\eta_n}_x(t),\alpha^{\eta_n}(t)) e^{-\lambda t} dt\\ 
            &\leq \lim_{n\rightarrow \infty} \left( \int_{0}^{T} 
            l( X^{\eta_n}_x(t),\alpha^{\eta_n}(t)) e^{-\lambda t} dt + 
            U_{A^{m}}^{\eta_n}( X^{\eta_n}_x(T))e^{-\lambda T} \right)\\
            &\leq \lim_{n\rightarrow \infty}  U_{A^m}^{\eta_n}(x) + 
            \lim_{n\rightarrow \infty}\frac{U_{A^m}^{+}(x)-  U_{A^m}^{\eta_n}(x)}{2}\\
            &\leq U_{A^{m}}^{+}(x).
    \end{aligned}
    \end{equation}

    \noindent\textbf{Step 3:} Taking limits in \eqref{inUmais11} as $T$ goes to infinity and using
    Lemma~\ref{TTTTT}, we obtain that 
    $$   
 	U_{A^{m}}^{+}(x)  \leq 
	\lim_{n\rightarrow \infty}  U_{A^{m}}^{\eta_n}(x)+ \lim_{n\to \infty}\frac{U_{A^{m}}^{+}(x)-  
	U_{A^{m}}^{\eta_n}(x)}{2} \leq U_{A^{m}}^{+}(x)
    $$
    which implies that 
    $$ \lim_{n\rightarrow \infty}  U_{A^{m}}^{\eta_n}(x)=U_{A^{m}}^{+}(x)\,,$$
    and the result follows by the monotonicity of $U^\eta_{A^m}$ with respect to $\eta$.
\end{proof}

Unfortunately, we are not able to prove that the previous convergence result holds for  unbounded control sets , but at
least we can prove that the $\eta$-value functions on compact control sets, $A^m$, converge to the value function $U^+_A$ on 
the unbounded control set, $A$.  

\begin{corollary}\label{etaaaaaaaaa}  
    Under the hypotheses in Lemma \ref{TTTTT}, there exists a sequence $(\eta_j,M_j)\to (0,\infty)$
    as $j\to\infty$, such that $$\lim_{j \to \infty}U_{A^{M_j}}^{\eta_j} (x)=U_A^{+}(x).$$
\end{corollary}

\begin{proof}
Given $ j \in \mathbb{N}$, there exists $m_j$ such that for $ m \geq m_j $, 
$$\vert U_{A}^{+}(x)- U_{A^{m}}^{+}(x)\vert \leq \frac{1}{2j}.$$
 Besides that,  by Theorem  \ref{cooo}, there exists $ \eta_j $ such that for  $ \eta'_j \leq \min\{\frac{1}{j},  \eta_j\}$  $$ \vert U_{A^{m_j}}^{+}(x) - U_{A^{m_j}}^{\eta'_j}(x) \vert \leq \frac{1}{2j}.$$
Consequently, 
$$\vert U_{A}^{+}(x)- U_{A^{m_j}}^{\eta_j}(x)\vert \leq  \vert U_A^{+}(x)-  U_{A^{m_j}}^{+}(x)\vert + \vert U_{A^{m_j}}^{+}(x) - U_{A^{m_j}}^{\eta_j}(x) \vert
\leq   \frac{1}{2j}+ \frac{1}{2j}= \frac{1}{j}. $$
Taking limits as  $j$ goes  to $\infty$ then $(\eta_j,m_j)\to (0,\infty)$, and we get the result. 
\end{proof}

%%%%%%%%%%%%%%%%%%%%%%%%%%%%%%%%%%%%%%%%%%%%%%%%%%%%%%%%%%%%
\section{Application to a superlinear example}
\label{examples} 
%%%%%%%%%%%%%%%%%%%%%%%%%%%%%%%%%%%%%%%%%%%%%%%%%%%%%%%%%%%%

The aim of this section is to give a quite general example of superlinear Hamiltonian and the
associated control problem satisfying the  hypotheses \HA, \HB, \HC, \LOCa, \LOCb, which
allow to apply all the results of this article.

We consider the following Hamiltonian: 
\begin{equation}\label{ex_non_lin}
	H_i(x,u,p) = \lambda u + d_i(x)^{r}|p|^{r}-f_i(x),
\end{equation}
where $r>1$, $d_i, f_i:\R^N\to \R$ are continuous functions in $\overline{\Omega}_i$ such that 
$$\begin{aligned}
    & 0<d_i(x)\leq |x|^{\kappa} \text{ is locally Lipschitz,}\\ 
    & 0 \leq f_i(x) \leq C_f|x|^{ a-\epsilon} + \overline{C}_f\,,\\
    & a\geq ( a-1+ \kappa) r\,,\quad \lambda\geq a^r\,.
\end{aligned}$$
Observe that, here, $\epsilon>0$ is a small constant. For instance, if the $d_i$ functions are bounded and $f_i$ are sublinear functions, we can take $\kappa=0$ and $a=1$ which
allows $r$ to be as big as we want, and $\lambda=1$.

%--------------------------------------
\subsection{The control problem}

This above Hamiltonian is obtained by considering the control problem deriving from the following dynamics
and costs: 
$$\begin{aligned}
    & b_i(x,\alpha_i):=c_{r} d_i(x)|\alpha_i|^{r-2}\alpha_i,,\\
    & l_i(x,\alpha_i):= f_i(x)+|\alpha_i|^{r}\,,\\
    & c_{r}=\frac{r}{\sqrt[r]{(r-1)^{r-1}}}\,.
\end{aligned}$$
 
 Observe that $b_i$ and $l_i$ satisfy hypotheses \HA, \HB and \HC. Indeed $b_i(x,0)=0$, and given
 $v$ in $\R^{N}$ different from zero and $ x \in \R^{N}$, there exists a control $\alpha_i$ such
 that $b_i(x,\alpha_i)=v$, where $\alpha_i$ is given by  $\alpha_i =\left(\frac{|v|}{c_{  r}
 d_i(x)}\right)^{\frac{1}{  r-1}} \frac{v}{|v|}$.  This implies that $\R^{N}\subset
 \mathcal{B}^{A_i}(x)$ for any $x\in \R^{N}$. Moreover, hypothesis \HC is also satisfied since 
 $$
 \lim_{|\alpha_i| \to \infty}\frac{l_i(x,\alpha_i)}{1+|b_i(x, \alpha_i)|}  
  =\lim_{|\alpha_i| \to \infty} \frac{f_i(x) +|\alpha_i|^{  r} }{c_{  r} d_i(x)|\alpha_i|^{  r-1}}=
  \infty\,,
$$
uniformly on compact subsets of $\R^N$.

%-----------------------------------
\subsection{Convexifying the problem}

Concerning the set $(b,l)$ we face a problem here. 
Though the sets $B_i(x)=\{b(x,\alpha_i):\alpha_i\in A_i\}$ and $L_i(x)$ (defined similarly) are
convex, nothing indicates that the set $BL_i(x):=\{(b_i(x,\alpha_i),l_i(x,\alpha_i))| \alpha_i \in
A_i\}$ is also convex.

In order to bypass this difficulty, we introduce the convex enveloppe of $BL$ and solve the
associated PDE problem. However, this approach requires some compactness argument to prove that by
enlarging $BL$ we get the same Hamiltonians.

So, let us consider the convex hull of $BL^{m}_i(x)$ where the exponent $m$ indicates that we are
taking the compact control sets $A^m:= B_m(0)\subset\R^N$. We now work with
$\mathcal{B}\mathcal{L}^{m}(x)$ defined as
\begin{equation}\label{tr}
 \mathcal{B}\mathcal{L}^{m}(x):= \left\{\begin{array}{ll} 
\overline{co}~(B^{m}_{1}(x)\times L^{m}_1(x)) &   \mbox{if } x_N>0,
\\ \overline{co}~(B^{m}_{2}(x)\times L^{m}_2(x)) & \mbox{if }  x_N < 0,
\\  \overline{co}((B^{m}_{1}(x)\times L^{m}_1(x)) \cup  (B^{m}_{2}(x)\times L^{m}_2(x)))   &  \mbox{if } x_N=0.
 \end{array} \right.
\end{equation}

We skip all the notations and details but this time, the control problem satisfies all the needed
hypotheses to get a suitable framework. We denote by $\bar H^m_i, \bar H^m_T$ the Hamiltonians
associated to the convex hull.

Now, in order to connect $\bar H^m_i, \bar H^m_T$ to $H^m_i$, $H^m_T$ we need the following result:
\begin{lemma}\label{Kconvexity} 
    Let $ F: \R^{N} \to \R$ be a convex function and $K \subset \R^{N}$ compact. Then
    $$\sup_{\Lambda \in \overline{co}(K)} F(\Lambda)= \sup_{ \Lambda \in K} F(\Lambda).$$
\end{lemma}

\begin{proof}
    For each $\epsilon >0$, we consider $F_{\epsilon}(\Lambda):= F(\Lambda)+ \epsilon|\Lambda|^{2}$
    which is strictly convex. Notice that $\overline{co}(K)$ is compact.
    The supremum of $F_{\epsilon}$ over $\overline{co}(K)$ is attained at a point
    $\Lambda_{\epsilon}$ in $\overline{co}(K)$.  From  Jensen’s inequality, for a strictly convex
    function,  we have
    $$\sup_{ \Lambda\in \overline{co}(K)} F_{\epsilon}(\Lambda) =F_{\epsilon}\left(\int_{K} \Lambda
    d\mu_{\epsilon} (\Lambda)\right)< \int_{K}F_{\epsilon}( \Lambda) d\mu_{\epsilon} (\Lambda))\leq
    \sup_{\Lambda \in K}F_{\epsilon}(\Lambda).$$ 
    So, necessarily the convex combination for $\Lambda_{\epsilon}$ is trivial, in other words
    $\mu_{\epsilon}=\delta_{\Lambda_{\epsilon}}$ and $\Lambda_{\epsilon} \in K$. Now,
    $$
        \sup_{\overline{co}(K)}F \leq  \sup_{\overline{co}(K)}F_{\epsilon}= F(\Lambda_{\epsilon})+
        \epsilon|\Lambda_{\epsilon}|^{2} \leq \sup_{K}F + \epsilon|\Lambda_{\epsilon}|^{2}
    $$ 
    and since K is compact, $\Lambda_{\epsilon} \to \Lambda_{*} \in K$ at least along a subsequence.
    So, passing to the limit as $\epsilon$ goes to $0$,  yields that
    $$ \sup_{\overline{co}(K)}F \leq \sup_{K}F. $$
    which gives the result. The other inequality is trivial.
\end{proof}

We apply this lemma to the Hamiltonians associated to compact control sets:
\begin{corollary}\label{cor6.2} The Hamiltonians $H_i^m, ~H^{m}_T$ and $H^{m,\mathrm{reg}}_T$ satisfy   
    $$H^{m}_i=\overline{H}^{m}_i,~~  H^{m}_T=\overline{H}^{m}_T,~~\mbox{ and }~~
    H^{m,\mathrm{reg}}_T=\overline{H}^{m,\mathrm{reg}}_T.$$
\end{corollary}
\begin{proof}
    Thanks to Lemma \ref{Kconvexity}, considering $F(b,l)= -b\cdot p -l $ and $K=BL^m_i(x)$,
    $K=BL^m_{T}(x)$ and $K=BL^{m, \mathrm{reg}}_{T}(x)$, respectively, we obtain the results.
\end{proof}

The last step consists in remembering that on each compact set $V\subset\R^N$, any supersolution of
the Ishii problem is actually a supersolution of $H^m_i$ for some $m\in\N_*$ while any stratified
subsolution is always a stratified subsolution of $H^m_i$ and $H^m_T$.

So, provided \LOCa and \LOCb are satisfied, the (LCR) is enough to ensure a comparison result. And
this (LCR) can be performed by using $H^m_i$, $H^m_T$ or equivalently, $\bar H^m_i$, $\bar H^m_T$.

This proves that the comparison result for stratified solutions works for the superlinear example.
And most associated results also follow.

Observe that thanks to Proposition \ref{H_i_cont} the supremum that defines $H_i$ and $H_T$ are attained 
on compact control sets under hypotheses \HA, \HB, \HC. From Corollary \ref{cor6.2} we know that the 
the Hamiltonians on compact control sets are equal to the Hamiltonias defined on the convex hull. 
Therefore, the problem is well defined using the convexifying arguments.

%---------------------------------------
\subsection{Checking \LOCa and \LOCb}

Checking \LOCa and \LOCb requires first to set some growth conditions. 
Let $\mathcal{C}$ be the set  of subsolutions of \eqref{ISHII},\eqref{ishiit} and supersolutions of
\eqref{ISHII}, $\omega$,  such that  for $\epsilon>0$
\begin{equation}\label{growthconditionnn}
|\omega(x)| \leq C_{\omega}|x|^{ a-\epsilon} + \overline{C}_{\omega}
\quad\mbox{with }~ a\geq ( a-1+ \kappa)r.  
\end{equation}
where $r,~\kappa $ are the constants that define the dynamic and cost functions. 
Lengthy but straightforward computations show that 
\begin{equation}
    \psi(x):=-(1+|x|^{2 a})^{1/2}
\end{equation}
is a stratified subsolution of \eqref{ISHII} for all $a>1/2$. Moreover, \LOCa is satisfied with
$u_\beta(x):=(1-\beta)u(x)-\beta\psi(x)$.
Concerning \LOCb, it is satisfied with 
 $$u_{\beta\gamma}(x) := u_{\beta}(x)- \gamma(1+|x-x_{\beta}|^{2})^{\frac{1}{2}}.
 $$ 

%-------------------------------------------------------------------
\subsection{Filippov approximations of the superlinear example.}

In this subsection we consider consider Filippov approximations. Let $\varphi:\mathbb{R}\rightarrow
[0,1]$ be a continuous functions satisfying $\lim\limits_{y\rightarrow \infty} \varphi(y)=1$ and
$\lim\limits_{y\rightarrow -\infty} \varphi(y)=0$. Let
$$\varphi_{\varepsilon}(y):=\varphi(\frac{y}{\varepsilon}).$$ 
We extend  $H_1$ and $H_2$ for any $x \in \R^{N}$ as follows

\begin{equation}\label{convervphi}
    H_{\varepsilon}(x,u,p):=\varphi_{\varepsilon}(x_N)H_1(x,u,p) +
    (1-\varphi_{\varepsilon}(x_N))H_2(x,u,p), ~~~ \forall x \in \mathbb{R}^{N}.  
\end{equation}

In the following result we prove that there exists a unique solution of the Filippov approximation
Hamiltonian, and this solution converges to the value function ${U}^-_A$. 

\begin{theorem}
    There exists a unique locally Lipschitz continuous solution $u_{\varepsilon}$ of
    \eqref{convervphi} satisfying 
    $$-C(1+|x|^{2  a})^{\frac{1}{2}-\epsilon} -C \leq u_{\varepsilon}(x)\leq 
    C(1+|x|^{2  a})^{\frac{1}{2}-\frac{\epsilon}{2 a}} +C$$
    where $C>\max\left\{C_m,  \frac{ C_f}{\lambda},  \frac{\overline{C}_f}{\lambda} \right\}$.
    Moreover, $u_\varepsilon$ converges to ${U}_{A}^{-}$ as  $\varepsilon $ goes to $0$,
    locally uniformly in~$\mathbb{R}^{N}$.  
\end{theorem}

\begin{proof}
    We only sketch the proof since it follows the arguments of \cite{manoel}: existence of solution
    for \eqref{convervphi} is obtained using the Perron Method (see \cite[p. 52]{guy}) in the
    following class of solutions:
    $$\begin{array}{ll}
        \mathcal{S}_{\varepsilon}:= \{u~\mbox{solution of \eqref{convervphi}} : 
    -C(1+|x|^{2  a})^{\frac{1}{2}-\epsilon} -C \leq u(x) \leq C(1+|x|^{2  a})^{\frac{1}{2}
    -\frac{\epsilon}{2 a}} + C~ \mbox{for all} ~ x \in \R^N\}.\end{array}  
    $$
    By standard comprison arguments, since $H_\epsilon$ is continuous and coercive, we know that
    there is at most one solution $u_\epsilon$.

    Finally, taking the limit of $u_{\varepsilon}$ when $\varepsilon$ goes to $0$ is done as in
    \cite[p. 31]{manoel}: this limit is a viscosity solution in the sense of Ishii of \eqref{ISHII}
    and moreover, since it satisfies the $H_T$ subsolution inequality, it is the unique stratified
    solution of \eqref{ISHII}. Therefore, by uniqueness the limit of $u_{\varepsilon}$ is $U_A^{-}$.
\end{proof}

%\bibliography{refer}

\end{document}